\documentclass{amsart}

\usepackage{amssymb}
\usepackage{epsfig}
\usepackage{color}

\usepackage[all]{xy}

\def\co{\colon\thinspace}
\def\R{\mathbb{R}}
\def\Z{\mathbb{Z}}

\def\N{\mathbb{N}}
\def\proofbox{\qed}

\newtheorem{theorem}{Theorem}[section]
\newtheorem{lemma}[theorem]{Lemma}
\newtheorem{corollary}[theorem]{Corollary}
\newtheorem{proposition}[theorem]{Proposition}
\newtheorem{claim}[theorem]{Claim}

\theoremstyle{definition}
\newtheorem{remark}[theorem]{Remark}
\newtheorem{definition}[theorem]{Definition}
\newtheorem{example}[theorem]{Example}
\newtheorem{question}[theorem]{Question}

\newcommand{\FA}{\ensuremath{(\mathrm{FA})}}
\newcommand{\QFA}{\ensuremath{(\mathrm{QFA})}}
\newcommand{\fa}[1]{\ensuremath{(\mathrm{FA}\mathcal{#1})}}
\newcommand{\qfa}[1]{\ensuremath{(\mathrm{QFA}\mathcal{#1})}}
\newcommand{\qfar}{\ensuremath{(\mathrm{QFA}\{\R\})}}

\title[Quasi-actions on trees and Property {\QFA}]{Quasi-actions on trees and Property {\QFA}}
\author{J. F. Manning, with Appendix by N. Monod and B. R\'emy}
\date{8 Dec 2005. Version 3.1.  MSC2000: 20F65 (primary), 20E08, 53C23, 22E40 (secondary)}

\begin{document}

\begin{abstract}
We prove some general results about quasi-actions on trees and define
Property \QFA, which is analogous to
Serre's Property \FA, but in the coarse setting.  
This property is shown to hold for a class of
groups, including $SL(n,\Z)$ for $n\geq 3$.  We also give a way of thinking
about Property \QFA\ by breaking it down into statements about particular
classes of trees.
\end{abstract} 

\maketitle

\tableofcontents

\section{Introduction}

Group quasi-actions 
are a natural coarse generalization of isometric group
actions (See Section \ref{section:definitions} for precise definitions.).
The main motivating question for this paper is:
\begin{question}\label{question:main}
What kind of finitely generated groups admit (or don't admit)
nontrivial quasi-actions on trees?
\end{question}  

Cobounded quasi-actions on bounded valence bushy trees were studied in
\cite{msw:quasiactI}, where it was shown that such quasi-actions are always
quasi-conjugate to isometric actions on trees.  
The same is not true for quasi-actions on
$\R$ or on infinite valence bushy trees.   Part of the reason for this is that
isometric actions on
$\R$-trees are always quasi-conjugate to actions on simplicial
trees, but this is not the complete story.  
Examples of quasi-actions on simplicial trees which are not
quasi-conjugate to actions on $\R$-trees are given in \cite{manning:cocycles}.
Such ``exotic'' quasi-actions on trees appear to be plentiful,  but it is not
clear how much information can be obtained from them.  We hope to
clarify the situation by offering some partial answers to Question
\ref{question:main}.  
  
Recall that a group $G$ is said to have 
Property \FA\ if for any isometric
action of $G$ on a simplicial 
tree $T$, there is some fixed point for the action (that is,
there is
some point $x\in T$ so that the orbit $Gx=\{x\}$).

\begin{definition}
We will say that a group $G$ has Property \QFA\ if for every 
quasi-action of $G$ on any tree $T$, there is some $x\in T$ so that the orbit 
$Gx$ has finite diameter (equivalently, every orbit has finite diameter).  
\end{definition}

Here is a brief outline.  Section \ref{section:definitions} consists mainly of
definitions and can probably be skipped by the expert.
In Section \ref{section:lemmas} we prove some useful facts about quasi-actions
on trees.  In Section \ref{section:main} we use these facts to prove Property
\QFA\ for a class of boundedly generated 
groups including $SL(n,\mathcal{O})$, for $n\geq 3$ and $\mathcal{O}$
the ring of integers of an algebraic number field.  In Section
\ref{section:coda} we try to understand Property \QFA\ by breaking it into
statements about different kinds of trees.  An Appendix by Nicolas
Monod and Bertrand R\'emy gives some examples of boundedly generated
lattices in Lie
groups which satisfy Kazhdan's Property $(T)$ but \emph{not} Property \QFA.  

Unless otherwise stated, all groups
are assumed to be finitely generated.

\subsection{Acknowledgements}
Thanks to my doctoral advisor Daryl Cooper for asking me (more than once!)
whether $SL(3,\Z)$
could quasi-act nontrivially on a tree.  
Thanks also to the Oxford
Mathematical Institute for hospitality while part of this work was being done.
Many people
have had helpful things to say to me about these problems, including
Nicolas Monod,  Kevin Whyte, and Hee Oh.   
This work was partially supported by an NSF Postdoctoral
Research Fellowship and by a UCSB Graduate Division Dissertation Fellowship.

 Some of these results have been
worked out independently but not published by Lee Mosher, Michah Sageev and
Kevin Whyte.

\section{Preliminaries}\label{section:definitions}

\subsection{Coarse geometry}
\begin{definition}
If $X$ and $Y$ are metric spaces, $K\geq 1$ and $C\geq 0$,
a \emph{$(K,C)$-quasi-isometric embedding} of $X$ into $Y$ is a function
$q\co X\to Y$ so that
 For all $x_1$, $x_2\in X$
\[\frac{1}{K}d(x_1,x_2)-C\leq d(q(x_1),q(x_2))\leq Kd(x_1,x_2)+C\]

If in addition the map $q$ is \emph{$C$-coarsely onto} -- \emph{i.e.},
every $y\in Y$ is
distance at most $C$ from some point in $q(X)$ -- then $q$ is called a
\emph{$(K,C)$-quasi-isometry}.
The two metric spaces $X$ and $Y$ are then
said to be \emph{quasi-isometric} to one
another.  This is a symmetric condition.
\end{definition}
\begin{definition}
A \emph{$(K,C)$-quasi-geodesic} in $X$
is a $(K,C)$-quasi-isometric embedding $\gamma\co\R\to X$.
We will occasionally abuse
notation by referring to the image of $\gamma$ as a quasi-geodesic.
\end{definition}
\begin{definition}
A \emph{$(K,C)$-quasi-action} of a group $G$ on a metric space $X$ is a map
$A\co G\times X\to X$, denoted $A(g,x)\mapsto gx$, so that
the following hold:
\begin{enumerate}
\item For each $g$,
$A(g,-)\co G\to G$ is a $(K,C)$-quasi-isometry.
\item For each $x\in X$ and $g$,
$h\in G$, we have \[d(g(hx), (gh)x)=d(A(g,A(h,x)),A(gh,x))\leq C.\]
\end{enumerate}
(Note that $K$ and $C$ must be independent of $g$ and $h$.)
 We call a quasi-action \emph{cobounded} if
for every $x\in X$, the map $A(-,x)\co G\to X$ is $C'$-coarsely onto
for some $C'$.
\end{definition}
\begin{definition}
Suppose that $A_X\co G\times X\to X$ and $A_Y\co G\times Y\to Y$ are
quasi-actions.  A map $f\co X\to Y$ is called coarsely equivariant if there is
some $C$ so that $d(f\circ A_X(g,x), A_Y(g,f(x)))\leq C$ for all $g$ in $G$ and
$x$ in $X$.

A coarsely equivariant quasi-isometry is called a \emph{quasi-conjugacy}.
\end{definition}

\begin{example}\label{example:quasicharacter}
Let $f\co G\to \R$ be a \emph{quasicharacter}; i.e., suppose that for
some $C\geq 0$ and for all $g$ and $h$ in $G$, $|\delta
f(g,h)|=|f(gh)-f(g)-f(h)|\leq C$ (see also Section
\ref{subsec:amenable}).
A $(1,C)$-quasi-action of $G$ on $\R$ is given by $A(g,x) = x+f(g)$.
\end{example}

\subsection{Quasi-trees and other hyperbolic spaces}
\begin{figure}[htbp]
\begin{center}
\begin{picture}(0,0)%
\includegraphics{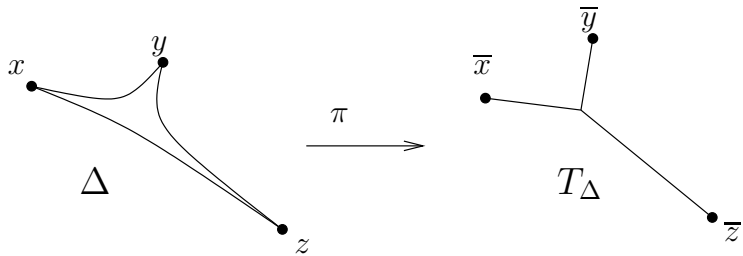}%
\end{picture}%
\setlength{\unitlength}{3947sp}%
\begingroup\makeatletter\ifx\SetFigFont\undefined%
\gdef\SetFigFont#1#2#3#4#5{%
  \reset@font\fontsize{#1}{#2pt}%
  \fontfamily{#3}\fontseries{#4}\fontshape{#5}%
  \selectfont}%
\fi\endgroup%
\begin{picture}(4500,1630)(451,-1910)
\put(451,-736){\makebox(0,0)[lb]{\smash{\SetFigFont{12}{14.4}{\familydefault}{\mddefault}{\updefault}{\color[rgb]{0,0,0}$x$}%
}}}
\put(1351,-586){\makebox(0,0)[lb]{\smash{\SetFigFont{12}{14.4}{\familydefault}{\mddefault}{\updefault}{\color[rgb]{0,0,0}$y$}%
}}}
\put(2251,-1861){\makebox(0,0)[lb]{\smash{\SetFigFont{12}{14.4}{\familydefault}{\mddefault}{\updefault}{\color[rgb]{0,0,0}$z$}%
}}}
\put(3376,-736){\makebox(0,0)[lb]{\smash{\SetFigFont{12}{14.4}{\familydefault}{\mddefault}{\updefault}{\color[rgb]{0,0,0}$\overline{x}$}%
}}}
\put(4051,-436){\makebox(0,0)[lb]{\smash{\SetFigFont{12}{14.4}{\familydefault}{\mddefault}{\updefault}{\color[rgb]{0,0,0}$\overline{y}$}%
}}}
\put(4951,-1786){\makebox(0,0)[lb]{\smash{\SetFigFont{12}{14.4}{\familydefault}{\mddefault}{\updefault}{\color[rgb]{0,0,0}$\overline{z}$}%
}}}
\put(2476,-1036){\makebox(0,0)[lb]{\smash{\SetFigFont{12}{14.4}{\familydefault}{\mddefault}{\updefault}{\color[rgb]{0,0,0}$\pi$}%
}}}
\put(901,-1486){\makebox(0,0)[lb]{\smash{\SetFigFont{14}{16.8}{\familydefault}{\mddefault}{\updefault}{\color[rgb]{0,0,0}$\Delta$}%
}}}
\put(3901,-1486){\makebox(0,0)[lb]{\smash{\SetFigFont{14}{16.8}{\familydefault}{\mddefault}{\updefault}{\color[rgb]{0,0,0}$T_\Delta$}%
}}}
\end{picture}
\caption{A triangle and its comparison tripod}
\label{fig:tripdef}
\end{center}
\end{figure}

All metric spaces will be assumed to be complete geodesic metric spaces, and the
distance between two points $x$ and $y$ will usually be denoted $d(x,y)$.
Several equivalent definitions and a much fuller discussion
of $\delta$-hyperbolic metric spaces can be
found, for instance, in \cite{bridhaef:book}, Chapter III.H.  We will use a
definition
which emphasizes the ``arboreality'' of hyperbolic spaces.  Given a geodesic
triangle $\Delta(x,y,z)$ in
any metric space, there is a unique \emph{comparison tripod}, $T_\Delta$,
a metric tree so
that the distances between the three extremal points of the tree,
$\overline{x}$, $\overline{y}$ and $\overline{z}$ , are the same as
the distances between $x$, $y$ and $z$ (See Figure \ref{fig:tripdef}.).
There is an obvious map $\pi\co \Delta(x,y,z)\to T_\Delta$ which takes $x$ to
$\overline{x}$, $y$ to $\overline{y}$ and $z$ to $\overline{z}$, and which is an isometry on
each side of $\Delta(x,y,z)$.
\begin{definition}
A space $X$ is \emph{$\delta$-hyperbolic} if for any geodesic triangle
$\Delta(x,y,z)$ and any point $p$ in the comparison tripod $T_\Delta$,
 the diameter of $\pi^{-1}(p)$ is less than $\delta$.  If $\delta$ is
 unimportant we may simply say that $X$ is \emph{Gromov hyperbolic}.
\end{definition}

\begin{definition}
Let $x$, $y$, $z\in X$.  The \emph{Gromov product} of $x$ and $y$ with 
respect to $z$
is $(x,y)_z=\frac{1}{2}(d(x,z)+d(y,z)-d(x,y))$.  Equivalently, $(x,y)_z$ is the
distance from $\overline{z}$ to the central vertex of the comparison tripod
$T_\Delta$ for any geodesic triangle $\Delta(x,y,z)$.
\end{definition}

The following is well known (see for example \cite{bridhaef:book},
III.H.1.22):
\begin{lemma}\label{lemma:products}
For any $\delta$ there is some $\delta'$ so that if $x$, $y$, $z$, and
$w$ are points in a $\delta$-hyperbolic space, then 
\[(x,y)_w\geq\min\{(x,z)_w,(y,z)_w\}-\delta'.\]
\end{lemma}

\begin{definition}
Fix some $z\in X$. 
We say that a sequence $\{x_i\}$
\emph{tends to infinity} if $\liminf_{i,j\to\infty}(x_i,x_j)_z =
\infty$.  On the set of such sequences we may define an equivalence
relation: $\{x_i\}\sim\{y_i\}$ if $\liminf_{i,j\to\infty}(x_i,y_i)_z =
\infty$.  
The \emph{Gromov boundary} of $X$, also written $\partial X$
is  the
set of equivalence classes of sequences tending to infinity.  The
Gromov boundary does not depend on the choice of $z$.  
\end{definition}
\begin{remark}
We may topologize $X\cup\partial X$ so that if $\{x_i\}$ tends to
infinity, then $\lim_{i\to\infty} x_i =[\{x_i\}]$.
Furthermore, if $\gamma\co{[0,\infty}\to X$ is a quasi-geodesic ray,
then for any sequence $\{t_i\}$ with $\lim_{i\to\infty}t_i=\infty$,
the sequence $\{\gamma(t_i)\}$ tends to infinity.  
The point $\{\gamma(t_i)\}\in\partial X$ does not depend on
the choice of $\{t_i\}$.  
\end{remark}

The following Lemma about ``stability'' of quasi-geodesics
is well-known.  We include a proof for
completeness and because we were unable to find this precise statement
in the literature (but see Remark \ref{remark:stability}).
\begin{lemma}\label{lemma:uniformnbhd}
Let $K\geq 1$, $C\geq 0$, $\delta\geq 0$.  Then there is some
$B=B(K,C,\delta)$, so that if $\gamma$ and $\gamma'$ are two
$(K,C)$-quasi-geodesics with the same endpoints in $X\cup\partial X$,
and $X$ is a $\delta$-hyperbolic geodesic metric space, then the image
of $\gamma$ lies in a $B$-neighborhood of the image of $\gamma'$.
\end{lemma}
\begin{proof}
Fix $K$, $C$, and $\delta$.  Several proofs of the lemma exist in the
literature under the assumption that $\gamma$ and $\gamma'$ are
quasi-geodesics of finite length (\cite{bridhaef:book}, III.H.1.7, for 
example).  Let $B_0=B_0(K,C,\delta)$ be the
constant which suffices in this case.  

Suppose first that $\gamma\co{[0,\infty)}\to X$ and
  $\gamma'\co{[0,\infty)}\to X$ are quasi-geodesic rays,
sharing one endpoint in $X$ and another in $\partial X$.  Let
$z=\gamma(0)=\gamma'(0)$.  
For each integer $i>0$ let $x_i$ be the first point in the
image of $\gamma$ which is distance $i$ from $z$, and let $y_i$ be
the first point in the image of $\gamma'$ which is a distance $i$ from
$z$.  We have $\{x_i\}\sim\{y_i\}$ and so in particular
$(x_i,y_i)_z\to\infty$ as $i\to\infty$.  

Note that there is an $\epsilon$ depending only on $K$, $C$, and
$\delta$ so that
the image of $\gamma$ is contained in an $\epsilon$-neighborhood of
$\{x_i\ |\ i\in\N\}$ and the image of $\gamma'$ is contained in an
$\epsilon$-neighborhood of $\{y_i\ |\ i\in\N\}$. 
(It is sufficient to bound $d(x_i,x_{i+1})$.
Let $p$ be a point on $[z,x_{i+1}]$ so
that $d(x_i,p)\leq B_0$.  Then:
\[d(x_i,x_{i+1})\leq d(p,x_i)+d(p,x_{i+1})\leq B_0+d(p,x_{i+1})\]
But since $d(p,z)+d(p,x_i)\geq d(x_i,z)=i$ we have
\[d(p,x_{i+1})= i+1-d(p,z) \leq i+1-(i-B_0) = B_0+1\]
which implies that $d(x_i,x_{i+1})\leq 2 B_0+1$.)  

Thus to prove the lemma, it suffices to bound $d(x_i, y_i)$ in
terms of $K$, $C$, and $\delta$.  Fixing $i$, choose $N$ so that
$(x_N,y_N)_z > i+B_0+2 \delta$ (See Figure \ref{fig:trick}.). 
\begin{figure}[htbp]
\begin{center}
\begin{picture}(0,0)%
\includegraphics{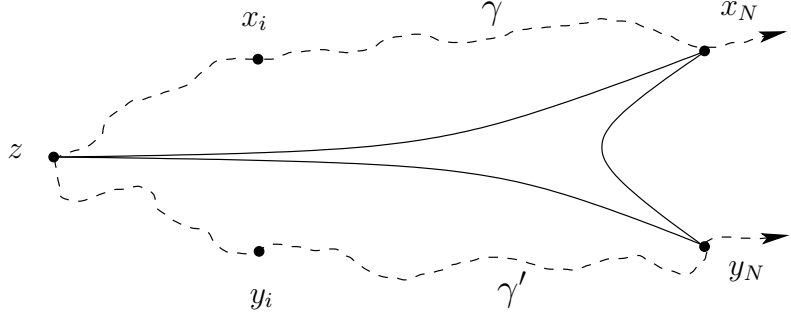}%
\end{picture}%
\setlength{\unitlength}{3947sp}%
\begingroup\makeatletter\ifx\SetFigFont\undefined%
\gdef\SetFigFont#1#2#3#4#5{%
  \reset@font\fontsize{#1}{#2pt}%
  \fontfamily{#3}\fontseries{#4}\fontshape{#5}%
  \selectfont}%
\fi\endgroup%
\begin{picture}(4962,2064)(871,-1790)
\put(886,-811){\makebox(0,0)[lb]{\smash{{\SetFigFont{12}{14.4}{\familydefault}{\mddefault}{\updefault}{\color[rgb]{0,0,0}$z$}%
}}}}
\put(3863, 67){\makebox(0,0)[lb]{\smash{{\SetFigFont{14}{16.8}{\familydefault}{\mddefault}{\updefault}{\color[rgb]{0,0,0}$\gamma$}%
}}}}
\put(3961,-1703){\makebox(0,0)[lb]{\smash{{\SetFigFont{14}{16.8}{\familydefault}{\mddefault}{\updefault}{\color[rgb]{0,0,0}$\gamma'$}%
}}}}
\put(2356, 22){\makebox(0,0)[lb]{\smash{{\SetFigFont{12}{14.4}{\familydefault}{\mddefault}{\updefault}{\color[rgb]{0,0,0}$x_i$}%
}}}}
\put(2408,-1718){\makebox(0,0)[lb]{\smash{{\SetFigFont{12}{14.4}{\familydefault}{\mddefault}{\updefault}{\color[rgb]{0,0,0}$y_i$}%
}}}}
\put(5363, 89){\makebox(0,0)[lb]{\smash{{\SetFigFont{12}{14.4}{\familydefault}{\mddefault}{\updefault}{\color[rgb]{0,0,0}$x_N$}%
}}}}
\put(5416,-1546){\makebox(0,0)[lb]{\smash{{\SetFigFont{12}{14.4}{\familydefault}{\mddefault}{\updefault}{\color[rgb]{0,0,0}$y_N$}%
}}}}
\end{picture}%
\caption{Closeness of quasi-geodesic rays follows from closeness of
long quasi-geodesic segments.}
\label{fig:trick}
\end{center}
\end{figure}
Let $p_x$ be some point on $[z,x_N]$ with $d(p_x,x_i)\leq B_0$, and
let $p_y$ be some point on $[z,y_N]$ with $d(p_y,y_i)\leq B_0$.  Note
that $d(z,p_x)$ and $d(z,p_y)$ are both at most $i+B_0$, and hence
less than $(x_N,y_N)_z$.  Consider the comparison tripod $T$
for a geodesic triangle with vertices $\{z,x_N,y_N\}$.
Both $\overline{p_x}$
and $\overline{p_y}$ lie in the leg of the tripod nearest
$\overline{z}$, as the length of this leg is precisely $(x_N,y_N)_z$.
(If $\xi$ is a point in the geodesic triangle, we write
$\overline{\xi}$ for the corresponding point in $T$.)
Thus 
\[d(\overline{p_x},\overline{p_y}) = |d(p_x,z)-d(p_y,z)| \leq 2 B_0.\]  
But this
implies that $d(p_x,p_y)\leq 2 B_0+\delta$, from which
it follows that $d(x_i,y_i)\leq 4 B_0 +\delta$.  It follows that any
point on $\gamma$ is within 
$B_1(K,C,\delta) = 4 B_0+\delta +\epsilon$ of the
$\gamma'$.  

Finally, suppose that both endpoints of $\gamma$ and $\gamma'$ are in
$\partial X$ and that
$\lim_{t\to\infty}\gamma(t)=\lim_{t\to\infty}\gamma'(t)$ and
$\lim_{t\to-\infty}\gamma(t)=\lim_{t\to-\infty}\gamma'(t)$. 
Let $D =
\inf_{s,t}d(\gamma(s),\gamma'(t))$ and re-parameterize $\gamma$ and
$\gamma'$ so that $d(\gamma(0),\gamma'(0)\leq D+1$.  It is not hard to
see that $\gamma$ must lie in a $B_1(K,C+D+1,\delta)$-neighborhood of
$\gamma'$.  Thus if we can find a universal bound $D_{\max}$
for $D$, we may set 
$B(K,C,\delta)=B_1(K,C+D_{\max}+1,\delta)$ and the lemma will be proved.

Since
$\liminf_{s,t\to\infty}(\gamma(s),\gamma'(t))_{\gamma(0)}=\infty$ we
may choose $s$ and $t$ so that
\[(\gamma(s),\gamma'(t))_{\gamma(0)}\geq D+1.\]
\begin{figure}[htbp]
\begin{center}
\begin{picture}(0,0)%
\includegraphics{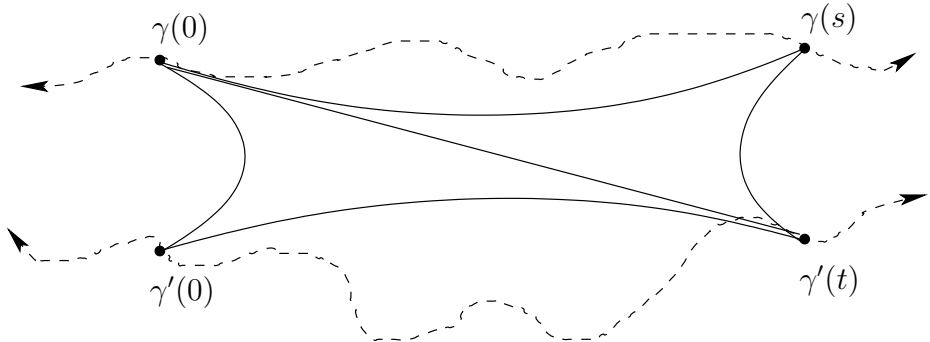}%
\end{picture}%
\setlength{\unitlength}{3947sp}%
\begingroup\makeatletter\ifx\SetFigFont\undefined%
\gdef\SetFigFont#1#2#3#4#5{%
  \reset@font\fontsize{#1}{#2pt}%
  \fontfamily{#3}\fontseries{#4}\fontshape{#5}%
  \selectfont}%
\fi\endgroup%
\begin{picture}(5844,2142)(394,-2015)
\put(1318,-1747){\makebox(0,0)[lb]{\smash{{\SetFigFont{12}{14.4}{\familydefault}{\mddefault}{\updefault}{\color[rgb]{0,0,0}$\gamma'(0)$}%
}}}}
\put(1331,-97){\makebox(0,0)[lb]{\smash{{\SetFigFont{12}{14.4}{\familydefault}{\mddefault}{\updefault}{\color[rgb]{0,0,0}$\gamma(0)$}%
}}}}
\put(5408,-32){\makebox(0,0)[lb]{\smash{{\SetFigFont{12}{14.4}{\familydefault}{\mddefault}{\updefault}{\color[rgb]{0,0,0}$\gamma(s)$}%
}}}}
\put(5397,-1679){\makebox(0,0)[lb]{\smash{{\SetFigFont{12}{14.4}{\familydefault}{\mddefault}{\updefault}{\color[rgb]{0,0,0}$\gamma'(t)$}%
}}}}
\end{picture}%
\caption{If $(\gamma(s),\gamma'(t))_{\gamma(0)}$ 
is large enough, the geodesics must come close together.}
\label{fig:two}
\end{center}
\end{figure}
In
this case any geodesic arcs between $\gamma(0)$ and $\gamma(s)$ and
between $\gamma'(0)$ and $\gamma'(t)$ must contain points which are at
most $2\delta$ apart (see Figure \ref{fig:two} and the comparison
tripods in Figure \ref{fig:two2}).
\begin{figure}[htbp]
\begin{center}
\begin{picture}(0,0)%
\includegraphics{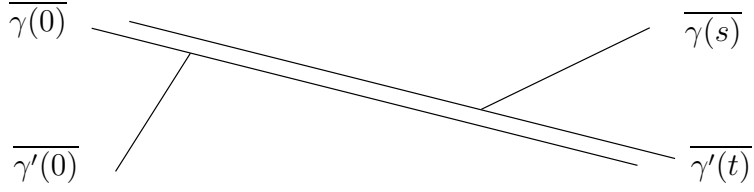}%
\end{picture}%
\setlength{\unitlength}{3947sp}%
\begingroup\makeatletter\ifx\SetFigFont\undefined%
\gdef\SetFigFont#1#2#3#4#5{%
  \reset@font\fontsize{#1}{#2pt}%
  \fontfamily{#3}\fontseries{#4}\fontshape{#5}%
  \selectfont}%
\fi\endgroup%
\begin{picture}(4309,1135)(1345,-1278)
\put(5639,-1208){\makebox(0,0)[lb]{\smash{{\SetFigFont{12}{14.4}{\familydefault}{\mddefault}{\updefault}{\color[rgb]{0,0,0}$\overline{\gamma'(t)}$}%
}}}}
\put(1360,-302){\makebox(0,0)[lb]{\smash{{\SetFigFont{12}{14.4}{\familydefault}{\mddefault}{\updefault}{\color[rgb]{0,0,0}$\overline{\gamma(0)}$}%
}}}}
\put(1386,-1209){\makebox(0,0)[lb]{\smash{{\SetFigFont{12}{14.4}{\familydefault}{\mddefault}{\updefault}{\color[rgb]{0,0,0}$\overline{\gamma'(0)}$}%
}}}}
\put(5604,-368){\makebox(0,0)[lb]{\smash{{\SetFigFont{12}{14.4}{\familydefault}{\mddefault}{\updefault}{\color[rgb]{0,0,0}$\overline{\gamma(s)}$}%
}}}}
\end{picture}%
\caption{The comparison tripods for the geodesic triangles in Figure
  \ref{fig:two} must fit together as shown, since
  $(\gamma(s),\gamma'(0))_{\gamma(0)}\leq D+1$, but
  $(\gamma(s),\gamma'(t))_{\gamma(0)}\geq D+1$.}
\label{fig:two2}
\end{center}
\end{figure}
There must therefore be points on the quasi-geodesic segments
$\gamma|_{[0,s]}$ and $\gamma'|_{[0,t]}$ which are within $2\delta +
B_0(K,C,\delta)$ of one another, and so $D_{\max}\leq
2\delta+B_0(K,C,\delta)$.  
\end{proof}
\begin{remark}\label{remark:stability}
A Gromov hyperbolic space $X$ is called \emph{ultra-complete} if every two points in
$X\cup\partial X$ are joined by a geodesic.  Trees are always
ultra-complete, but general locally infinite Gromov hyperbolic graphs
need not be.  It is claimed
in Section 7.5 of \cite{gromov:wordhyperbolic} that any
$\delta$-hyperbolic geodesic metric space $X$ isometrically embeds in an
ultra-complete space $Y$ with $\sup_{y\in Y}d(y,X)\leq C$ for some 
$C<\infty$.  
Another way to prove Lemma \ref{lemma:uniformnbhd} would be to show first
that $C$ depends only on $\delta$ and then apply 
Th\'eor\`eme 3.1 of \cite{cdp}. 
\end{remark}

\begin{definition}\label{definition:ellhyp}
Fix $x\in X$, where $X$ is a $\delta$-hyperbolic metric space on which
$G$ quasi-acts.
Let $O_{g,x}\co\R\to X$ be defined by $O_{g,x}(t)=g^{\lfloor
  t\rfloor}x$,
where $\lfloor t\rfloor$ is
the largest integer smaller than  $t$.
If $O_{g,x}$ has
bounded image, we say $g$ quasi-acts \emph{elliptically}.  If $O_{g,x}$ is a
quasi-geodesic, then we say $g$ quasi-acts \emph{hyperbolically}.  (If
$G$ acts isometrically on $X$, we may simply say $g$ \emph{acts}
elliptically or hyperbolically.)
\end{definition}
It is not hard to check that Definition \ref{definition:ellhyp} is
 independent of $x$ and
agrees with the standard definitions in case $G$ acts isometrically.
It has the added benefit of being invariant under quasi-conjugacy.  In the case
that $G$ acts isometrically on $X$, and $g\in G$ acts hyperbolically, then $g$
always has a \emph{quasi-axis}, a quasi-geodesic whose image is invariant under
the infinite cyclic group $\langle g \rangle$.  Indeed if $x\in X$, and
$\gamma_0\co[0,1]\to X$ is a geodesic segment with $\gamma_0(0)=x$ and
$\gamma_0(1)=gx$, then the reader can easily verify that $\gamma\co\R\to X$ is a
continuous quasi-geodesic, if we define
$\gamma(t)=g^{\lfloor t\rfloor}\gamma_0(t-\lfloor t\rfloor)$.

\begin{example}\label{example:quasifrompseudo}
If the quasicharacter
$f\co G\to\R$ from Example \ref{example:quasicharacter} above is a homomorphism
restricted to each cyclic subgroup, it is called a \emph{pseudocharacter} or
\emph{homogeneous} quasicharacter.  As in Example \ref{example:quasicharacter},
the pseudocharacter $f$ induces a quasi-action of $G$ on $\R$.
A group element $g\in G$ quasi-acts hyperbolically if
and only if $f(g)$ is nonzero.  Unless $f$ is identically $0$, the quasi-action
is cobounded.
\end{example}

\begin{definition}
A \emph{quasi-tree} is a complete geodesic
metric space quasi-isometric to some simplicial tree
(All simplicial trees are assumed to be endowed with a path metric in which
every edge has length 1.).
\end{definition}

Quasi-trees satisfy a particularly strong form of $\delta$-hyperbolicity.
\begin{lemma}\label{lemma:bottleneck}
If $X$ is a quasi-tree, then there is a $\delta>0$ so that:
\begin{enumerate}
\item\label{bottle}For
any two points $x$ and $y$ in $X$, and any point $p$ on a geodesic between
$x$ and $y$, any path from $x$ to $y$ must pass within $\delta$ of $p$.
\item\label{hyper} $X$ is $\delta$-hyperbolic.
\end{enumerate}
\end{lemma}
\begin{proof}
Exercise.
\end{proof}
\begin{definition}
If a quasi-tree $X$ satisfies the conclusions of Lemma
\ref{lemma:bottleneck} for $\delta\geq 0$, we say that $X$ is a
$\delta$-quasi-tree.
\end{definition}

\subsection{Bounded cohomology, amenability, and Trauber's
Theorem}\label{subsec:amenable}
We give only a few needed facts here.  For fuller
discussion of these topics, see \cite{gromov:bounded} and
\cite{grigorchuk:bounded}.

We give the definition of bounded cohomology for groups only.
\begin{definition}
The \emph{bounded cohomology} $H_b^*(G;\R)$ of a group $G$
is the cohomology of the cochain complex
$C_b^*(G;\R)$, where
\[C_b^n(G;\R)=\{f\co G^n\to \R\ |\ \sup_{G^n}|f(g_1,\ldots,g_n)|<\infty\}\]
 and $\delta\co C_b^n(G;\R)\to C_b^{n+1}(G;\R)$ is given by
\begin{eqnarray*}
\delta f(g_1,\ldots,g_{n+1}) & = & f(g_2,\ldots,g_{n+1})
+\Sigma_{i=1}^n(-1)^i f(g_1,\ldots
g_{i-1},g_i g_{i+1},\ldots,g_{n+1})\\
& & \mbox{} + (-1)^{n+1} f(g_1,\ldots,g_n)
\end{eqnarray*}
This cochain complex is a sub-complex of the complex 
$C^*(G;\R)$ of \emph{all}
real valued functions on $G$, $G\times G$, and so on.  The cohomology of
$C^*(G;\R)$ is the ordinary cohomology of $G$ with real coefficients.
\end{definition}

\begin{definition}
A \emph{quasicharacter} is an element $f$ of $C^1(G;\R)$ whose coboundary
$\delta f$ lies in $C^2_b(G;\R)$.
The quasicharacter $f$ is a
\emph{pseudocharacter} if in addition $f(g^n)=nf(g)$ for all $n\in \Z$ and $g\in
G$.  In either case we define the \emph{defect} of $f$ as
$\|\delta f\|=\sup_{\mbox{$g$, $h\in G$}}|\delta
f(g,h)|$.
\end{definition}

\begin{remark}\label{remark:pseudoquasi}
Note that if $f$ is a quasicharacter, and $\phi$ is given by
$\phi(g)=\lim_{n\to\infty}\frac{f(g^n)}{n}$,
then $\phi$ is a pseudocharacter with
$\phi-f$ bounded and $[\delta \phi]=[\delta f]$.
If a pseudocharacter $\phi$ is ever nonzero, then
it is unbounded.
\end{remark}
The relationship between quasicharacters and pseudocharacters and bounded
cohomology is a major tool for understanding $H^2_b$ in certain situations (see
for example
 \cite{bestvinafujiwara:mcg} and
\cite{grigorchuk:bounded}).

We will need only a few facts about amenable groups.
First, nilpotent groups are
amenable.  Second, amenable groups contain no free subgroups. Third:
\begin{theorem}\label{theorem:trauber}
(Trauber's Theorem) If $G$ is an amenable group then $H_b^n(G;\R)=0$ for all
$n$.
\end{theorem}
For a definition of amenability and a proof of Theorem
\ref{theorem:trauber} see \cite{grigorchuk:bounded}.

\section{Lemmata}\label{section:lemmas}

This section contains some general results about quasi-actions on trees
by (finitely generated) groups.  The key idea is that quasi-actions on 
trees and isometric actions on quasi-trees are essentially equivalent.  
  Proposition \ref{prop:quasitree}
gives a way to replace a quasi-action on a tree by an isometric action on a
Cayley graph which is a quasi-tree.  In \ref{subsec:Z} it is shown that there is
no such thing as a ``parabolic'' isometry of a quasi-tree.  In
\ref{subsec:pseudo} we show how to obtain a pseudocharacter from a quasi-action
on a tree which fixes one end.  

\subsection{Getting some action}\label{subsec:action}
Recall that if $G$ is a group and $S$ some (not necessarily finite) 
generating set, then we may form the Cayley graph $\Gamma(G,S)$ by setting
the zero-skeleton
$\Gamma(G,S)^0=G$ and connecting $g$ to $gs$ with an edge whenever $s\in S$.
We make $\Gamma(G,S)$ a metric space with a path metric in which every edge has
length $1$.  Then $G$ acts on the left by isometries of $\Gamma(G,S)$.

\begin{proposition}\label{prop:quasitree}
Suppose a finitely generated group $G$ quasi-acts  on
 a simplicial
tree $T$.  Then there is a generating set $S$
 for $G$ so that 
the Cayley graph $\Gamma(G,S)$
 embeds coarsely equivariantly and
quasi-isometrically in $T$.  Specifically, for $x\in T$ and $R$ 
sufficiently
large, we may take $S=\{s\in G\ |\  d(s(x),x)\leq R\}$.
\end{proposition}

\begin{proof}
We suppose that $G$ has a $(K,C)$-quasi-action on the simplicial tree $T$.  Let
$S_0$ be a finite generating set for $G$, and let $\Gamma_0=\Gamma(G,S_0)$ be
the associated Cayley graph.  Fix $x\in T$ and 
define $\pi_0\co \Gamma_0\to T$ by $\pi_0(g)=gx$.  We may assume $\pi_0$ sends
each edge to a geodesic.  For $s\in S_0$ and $g\in G$ we have:
\[d(gx,(gs)x)\leq d(gx,g(sx))+C\leq K d(x,sx)+2C\]
Let $D =  K(\sup_{s\in S_0}  d(x,sx))+2C$.  We fix  $R\geq 2KD+KC$.

We now set $S=\{s\in G\ |\ d(s(x),x)\leq R\}$, as in the
statement of the lemma.  We let $\Gamma=\Gamma(G,S)$ and extend $\pi_0$ to
$\pi\co\Gamma\to T$ which we may also assume sends each edge to a geodesic.  
To show that $\pi$ is a quasi-isometric
embedding, we must bound $d(\pi(p),\pi(q))$ above and below by affine functions
of $d(p,q)$.  We may restrict our attention to the case when both $p$ and $q$
are vertices (group elements), as $\Gamma$ is quasi-isometric to its
zero-skeleton.  

Suppose that $d(p,q)=1$.  Then there is some $s\in S$ so that
$p=qs$ and so 
\[d(px,qx)=d(qsx,qx)\leq d(q(sx),qx)+C \leq K d(sx,x)+2C \leq KR+2C.\]
In general we have $d(\pi(p),\pi(q))=d(px,qx)\leq (KR+2C) d(p,q)$.  

In the other direction, note first that if $d(\pi(p),\pi(q))<R/K-C$, 
then $d(p,q)=1$.  If on the other hand
$d(\pi(p),\pi(q))\geq R/K-C$, we will choose a new $p'$ so that $d(p,p')=1$ and
$\pi(p')$ is closer to $\pi(q)$ than $\pi(p)$ was.
Note that $d(\pi(p),\pi(q))\geq R/K-C$ implies
$d(\pi(p),\pi(q))\geq 2D$ by our choice of $R$.
Let $z$ be the point on the
geodesic $[\pi(p),\pi(q)]$ which is a distance of $3D/2$ from $\pi(p)$.  This
point is in the image of $\Gamma(G,S)$, and so there is a group element $p'$ so
that $d(\pi(p'),z)\leq D/2$.  Since $d(p'x,p x)\leq 2D$ we have
$d(p^{-1}p'x,x)\leq 2KD+C<R$, implying $p^{-1}p'\in S$ and so $d(p',p)=1$.  On
the other hand, $\pi(p')$ is at least $D$ closer to $\pi(q)$ than $\pi(p)$ is.
Thus we can travel from $p$ to $q$ in $\Gamma$ by traversing at most
$d(px,qx)/D+1$ edges.  In other words, $d(p,q)\leq d(px,qx)/D + 1$ or
$d(\pi(p),\pi(q))=d(px,qx)\geq D d(p,q)-D$.  

Since 
\[ D d(p,q)-D\leq d(\pi(p),\pi(q))\leq (KR+2C)d(p,q),\]
the map $\pi$ is a quasi-isometric embedding.

We now show $\pi$ is coarsely equivariant.  
Again, we may restrict attention to vertices of
$X$.  Let $p$ be a vertex of $X$, and let $g\in G$.  We need a universal bound
on $d(g(\pi(p)), \pi(g(p))$.  Since $p$ is a vertex of $X$ it is a group
element, and so $\pi(g(p))=\pi(gp) = (gp)x$ and $\pi(p)=px$.  By the definition
of a $(K,C)$-quasi-action, $d(g(\pi(p)),\pi(g(p)))=d(g(px),(gp)x)\leq C$, and so
$\pi$ is coarsely equivariant.  
\end{proof}
\begin{remark}
Note that the isometric action we obtain from \ref{prop:quasitree} is
quasi-conjugate to the original quasi-action only in the cobounded
case 
(compare
Proposition 4.4 of \cite{manning:cocycles}).
This disadvantage is balanced by the fact that we may now work with a
left-invariant metric on $G$ itself.  Whether or not the original
quasi-action is cobounded, the Cayley graph $\Gamma(G,S)$ is
quasi-isometric to its image in $T$ and is hence a quasi-tree.
Conversely, if $G$ acts by isometries on a quasi-tree, then there is a
quasi-conjugate quasi-action on a tree.
\end{remark}

\subsection{Cayley graphs of $\Z$}\label{subsec:Z}
The aim in this subsection is to explain why every element of a group
quasi-acting on a tree must quasi-act either hyperbolically or
elliptically, in the sense of Definition \ref{definition:ellhyp}.
Suppose $g\in G$ and that $G$ quasi-acts on the tree $T$.  Then the
integers also quasi-act on $T$, via $n (x) = g^n x$ for $n\in\Z$.  By
Proposition \ref{prop:quasitree}, there is therefore a coarsely
equivariant quasi-isometric embedding of some Cayley graph 
$\Gamma=\Gamma(\Z,S)$ into $T$.  In particular this Cayley graph must
be $\delta$-hyperbolic for some $\delta$.  The element $g$ quasi-acts
elliptically if and only if the diameter of $\Gamma$ is finite, and
quasi-acts hyperbolically if and only if $S$ is finite.  The following
proposition shows that there are no other possibilities.
\begin{proposition}\label{prop:Z}
Suppose that $\Gamma=\Gamma(\Z,S)$ is $\delta$-hyperbolic.  Then
either $S$ is finite or the diameter of $\Gamma$ is finite.  
\end{proposition}
\begin{proof}
Let $\Delta\geq \max\{\delta, \delta', 1\}$, where $\delta'$ is the
constant from Lemma \ref{lemma:products}.
We will write $|n|_S$ for the distance $d(0,n)$ in $\Gamma$.
Assuming that the diameter if $\Gamma$ is infinite, we fix some $N>0$ so
that $|N|_S\geq 10\Delta$.  For $k\in\Z$, let $D_k=|kN|_S$. Note that
$D_{-k}=D_k$ for any $k$.     
\begin{claim}
$D_k\geq |k|(|N|_S - 4\Delta)$.
\end{claim}
\begin{proof}
It suffices to prove the Claim for all positive $k$.  We inductively
argue that the following two assertions hold for each $k$:
\begin{description}
\item [$P_k$:] \label{long} $D_k\geq k(|N|_S-4\Delta)$.
\item [$Q_k$:] \label{short} $(-N\ ,\ (kN))_0\leq 2\Delta$.
\end{description}
The statement $P_1$ is obvious; the statement $Q_1$ can be proved as
follows:
Let $\sigma$ be a geodesic segment from $-N$ to $0$, and let $T$ be a
geodesic triangle two of whose sides are $\sigma$ and $\sigma+N$.  
Let $z$ be a point on $\sigma$ so that 
$d(z,0)=(-N\ ,\ N)_0$ (see Figure \ref{figure:basecase}).  
\begin{figure}[htbp]
\begin{center}
\begin{picture}(0,0)%
\includegraphics{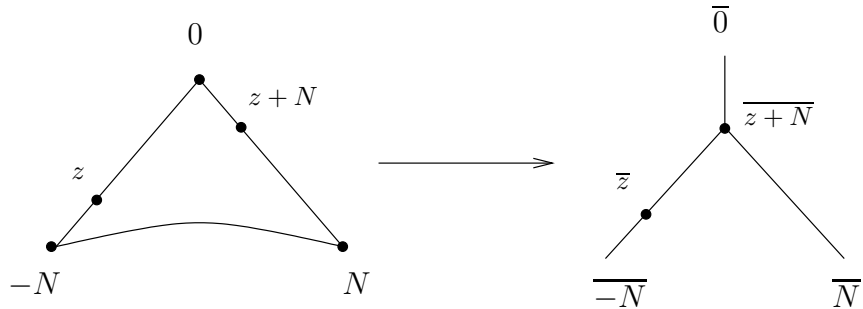}%
\end{picture}%
\setlength{\unitlength}{3947sp}%
\begingroup\makeatletter\ifx\SetFigFont\undefined%
\gdef\SetFigFont#1#2#3#4#5{%
  \reset@font\fontsize{#1}{#2pt}%
  \fontfamily{#3}\fontseries{#4}\fontshape{#5}%
  \selectfont}%
\fi\endgroup%
\begin{picture}(5277,1969)(361,-1709)
\put(4801, 89){\makebox(0,0)[lb]{\smash{{\SetFigFont{12}{14.4}{\familydefault}{\mddefault}{\updefault}{\color[rgb]{0,0,0}$\overline{0}$}%
}}}}
\put(4051,-1636){\makebox(0,0)[lb]{\smash{{\SetFigFont{12}{14.4}{\familydefault}{\mddefault}{\updefault}{\color[rgb]{0,0,0}$\overline{-N}$}%
}}}}
\put(5551,-1636){\makebox(0,0)[lb]{\smash{{\SetFigFont{12}{14.4}{\familydefault}{\mddefault}{\updefault}{\color[rgb]{0,0,0}$\overline{N}$}%
}}}}
\put(4193,-901){\makebox(0,0)[lb]{\smash{{\SetFigFont{10}{12.0}{\familydefault}{\mddefault}{\updefault}{\color[rgb]{0,0,0}$\overline{z}$}%
}}}}
\put(4996,-496){\makebox(0,0)[lb]{\smash{{\SetFigFont{10}{12.0}{\familydefault}{\mddefault}{\updefault}{\color[rgb]{0,0,0}$\overline{z+N}$}%
}}}}
\put(1501, 14){\makebox(0,0)[lb]{\smash{{\SetFigFont{12}{14.4}{\familydefault}{\mddefault}{\updefault}{\color[rgb]{0,0,0}$0$}%
}}}}
\put(376,-1561){\makebox(0,0)[lb]{\smash{{\SetFigFont{12}{14.4}{\familydefault}{\mddefault}{\updefault}{\color[rgb]{0,0,0}$-N$}%
}}}}
\put(2476,-1561){\makebox(0,0)[lb]{\smash{{\SetFigFont{12}{14.4}{\familydefault}{\mddefault}{\updefault}{\color[rgb]{0,0,0}$N$}%
}}}}
\put(773,-833){\makebox(0,0)[lb]{\smash{{\SetFigFont{10}{12.0}{\familydefault}{\mddefault}{\updefault}{\color[rgb]{0,0,0}$z$}%
}}}}
\put(1883,-368){\makebox(0,0)[lb]{\smash{{\SetFigFont{10}{12.0}{\familydefault}{\mddefault}{\updefault}{\color[rgb]{0,0,0}$z+N$}%
}}}}
\end{picture}%
\caption{Bounding $(-N,\ N)_0$.}
\label{figure:basecase}
\end{center}
\end{figure}
Then $d(z,z+N)$ is at most
$|N|-2(-N\ ,\ N)_0+\Delta$. 
On the other hand, since $N$ moves every vertex of $\Gamma(\Z,S)$ the
same distance, we must have $d(z,z+N)\geq |N|-1$ ($z$ may lie in the
middle of an edge).  Combining these two
inequalities gives us $2(-N\ ,\ N)_0\leq \Delta+1$, and so
\[(-N\ ,\ N)_0\leq\Delta\leq 2\Delta\] 
(assertion $Q_1$).  

Assuming $k\geq 2$, we argue the inductive step as follows:
Note first that 
$D_k=d(-N,(k-1)N) = |N|_S+|(k-1)N|_S-2(-N,\ (k-1)N)_0$.  Since
$(-N,(k-1)N)\leq 2\Delta$ by the induction hypothesis $Q_{k-1}$, we
have $D_k\geq D_{k-1}+|N|_S-4\Delta$.  Applying $P_{k-1}$ yields
$D_k\geq k(|N|_S-4\Delta)$, and the assertion $P_k$ is proved.  

To prove the statement $Q_k$, consider four geodesic triangles as pictured in
Figure \ref{figure:induction}
\begin{figure}[htbp]
\begin{center}
\begin{picture}(0,0)%
\includegraphics{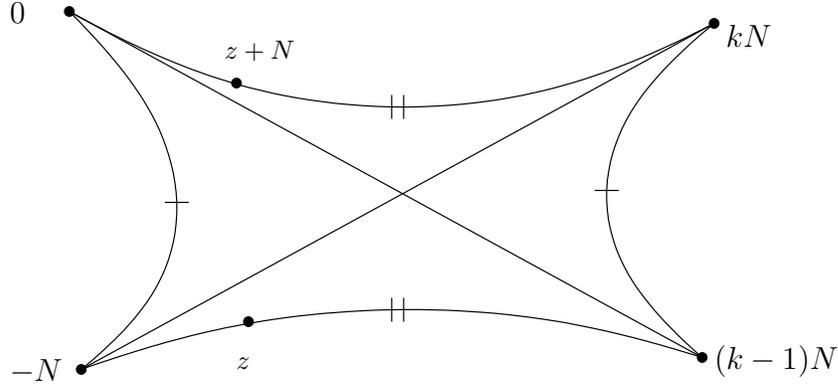}%
\end{picture}%
\setlength{\unitlength}{3947sp}%
\begingroup\makeatletter\ifx\SetFigFont\undefined%
\gdef\SetFigFont#1#2#3#4#5{%
  \reset@font\fontsize{#1}{#2pt}%
  \fontfamily{#3}\fontseries{#4}\fontshape{#5}%
  \selectfont}%
\fi\endgroup%
\begin{picture}(4530,2478)(-14,-3055)
\put(  1,-736){\makebox(0,0)[lb]{\smash{{\SetFigFont{12}{14.4}{\familydefault}{\mddefault}{\updefault}{\color[rgb]{0,0,0}$0$}%
}}}}
\put(  1,-2986){\makebox(0,0)[lb]{\smash{{\SetFigFont{12}{14.4}{\familydefault}{\mddefault}{\updefault}{\color[rgb]{0,0,0}$-N$}%
}}}}
\put(4426,-2911){\makebox(0,0)[lb]{\smash{{\SetFigFont{12}{14.4}{\familydefault}{\mddefault}{\updefault}{\color[rgb]{0,0,0}$(k-1)N$}%
}}}}
\put(4501,-886){\makebox(0,0)[lb]{\smash{{\SetFigFont{12}{14.4}{\familydefault}{\mddefault}{\updefault}{\color[rgb]{0,0,0}$kN$}%
}}}}
\put(1426,-2911){\makebox(0,0)[lb]{\smash{{\SetFigFont{10}{12.0}{\familydefault}{\mddefault}{\updefault}{\color[rgb]{0,0,0}$z$}%
}}}}
\put(1351,-961){\makebox(0,0)[lb]{\smash{{\SetFigFont{10}{12.0}{\familydefault}{\mddefault}{\updefault}{\color[rgb]{0,0,0}$z+N$}%
}}}}
\end{picture}%
\caption{Induction step.}
\label{figure:induction}
\end{center}
\end{figure}
and the corresponding comparison tripods shown in Figure
\ref{figure:peds}.  The geodesic segment between $-N$ and $(k-1)N$ 
in Figure \ref{figure:induction} should be 
chosen to be the obvious
translate of the one between $0$ and $kN$.   
\begin{figure}[htbp]
\begin{center}
\begin{picture}(0,0)%
\includegraphics{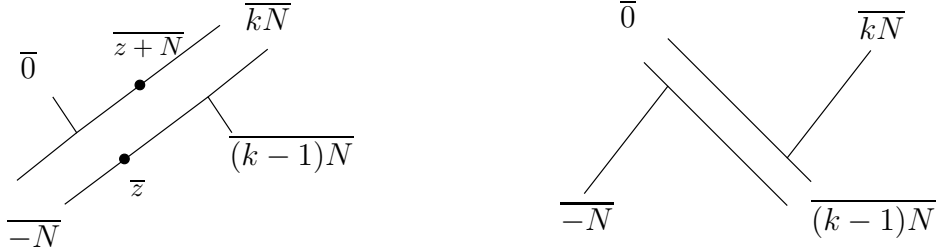}%
\end{picture}%
\setlength{\unitlength}{3947sp}%
\begingroup\makeatletter\ifx\SetFigFont\undefined%
\gdef\SetFigFont#1#2#3#4#5{%
  \reset@font\fontsize{#1}{#2pt}%
  \fontfamily{#3}\fontseries{#4}\fontshape{#5}%
  \selectfont}%
\fi\endgroup%
\begin{picture}(5450,1631)(188,-2496)
\put(5251,-2311){\makebox(0,0)[lb]{\smash{{\SetFigFont{12}{14.4}{\familydefault}{\mddefault}{\updefault}{\color[rgb]{0,0,0}$\overline{(k-1)N}$}%
}}}}
\put(3676,-2311){\makebox(0,0)[lb]{\smash{{\SetFigFont{12}{14.4}{\familydefault}{\mddefault}{\updefault}{\color[rgb]{0,0,0}$\overline{-N}$}%
}}}}
\put(4051,-1036){\makebox(0,0)[lb]{\smash{{\SetFigFont{12}{14.4}{\familydefault}{\mddefault}{\updefault}{\color[rgb]{0,0,0}$\overline{0}$}%
}}}}
\put(5551,-1111){\makebox(0,0)[lb]{\smash{{\SetFigFont{12}{14.4}{\familydefault}{\mddefault}{\updefault}{\color[rgb]{0,0,0}$\overline{kN}$}%
}}}}
\put(203,-2423){\makebox(0,0)[lb]{\smash{{\SetFigFont{12}{14.4}{\familydefault}{\mddefault}{\updefault}{\color[rgb]{0,0,0}$\overline{-N}$}%
}}}}
\put(1598,-1906){\makebox(0,0)[lb]{\smash{{\SetFigFont{12}{14.4}{\familydefault}{\mddefault}{\updefault}{\color[rgb]{0,0,0}$\overline{(k-1)N}$}%
}}}}
\put(1696,-1059){\makebox(0,0)[lb]{\smash{{\SetFigFont{12}{14.4}{\familydefault}{\mddefault}{\updefault}{\color[rgb]{0,0,0}$\overline{kN}$}%
}}}}
\put(286,-1358){\makebox(0,0)[lb]{\smash{{\SetFigFont{12}{14.4}{\familydefault}{\mddefault}{\updefault}{\color[rgb]{0,0,0}$\overline{0}$}%
}}}}
\put(976,-2101){\makebox(0,0)[lb]{\smash{{\SetFigFont{10}{12.0}{\familydefault}{\mddefault}{\updefault}{\color[rgb]{0,0,0}$\overline{z}$}%
}}}}
\put(870,-1216){\makebox(0,0)[lb]{\smash{{\SetFigFont{10}{12.0}{\familydefault}{\mddefault}{\updefault}{\color[rgb]{0,0,0}$\overline{z+N}$}%
}}}}
\end{picture}%
\caption{Comparison tripods for Figure \ref{figure:induction}.  The
  tripods on the left are congruent to one another, as are those on
  the right.}
\label{figure:peds}
\end{center}
\end{figure} 
Lemma \ref{lemma:products} implies that 
\begin{equation}\label{eq:prod}
(-N,\ (k-1)N)_0+\Delta\geq \min\{(-N,\ kN)_0,(kN,\  (k-1)N)_0\}.
\end{equation}
By assertion $Q_{k-1}$, the left hand side of equation \ref{eq:prod}
can be no larger than $3\Delta$. 
It follows from the congruence of the two right hand tripods in figure
\ref{figure:peds} that 
$(kN,\ (k-1)N)_0 = D_{k-1}-(kN,\ 0)_{(k-1)N} = 
D_{k-1}-(-N,\ (k-1)N)_0$.  
Thus by $Q_{k-1}$ and $P_{k-1}$, 
$(kN,\ (k-1)N)_0$ is at least $4\Delta$ and so equation \ref{eq:prod}
implies that 
$(-N,\ kN)_0\leq 3\Delta$.  This doesn't yet establish $Q_k$, but it
at least shows that the comparison tripods in the left half of Figure
\ref{figure:peds} are qualitatively correct;
since $|N|_S\geq 2 (-N,\ kN)_0$ and $kN > |N|_S$,
the lengths of the three legs of the leftmost tripod
must be ordered:
\[(-N,\ kN)_0\leq (0,\ kN)_{-N}\leq(0,\ -N)_{kN}\]
Now let $z$ be the point on $[-N,(k-1)N]$ so that 
$d(z,-N)=(0.\ kN)_{-N}$.  Then we have 
\begin{equation}\label{eq:quadz}
|N|_S-1\leq d(z,z+N)\leq 2\Delta+(0,kN)_{-N} - (-N,kN)_0.
\end{equation}
Since $(0,\ kN)_{-N} = |N|_S-(-N,\ kN)_0$, equation \ref{eq:quadz}
may be rewritten as 
\begin{equation}
|N|_S-1\leq 2\Delta+|N|_S - 2(-N,\ kN)_0,
\end{equation}
which can be rearranged as $(-N,\ kN)_0\leq \Delta + \frac{1}{2}\leq
2\Delta$, and $Q_k$ is proved.
\end{proof}

We now show that $S$ is finite.
Suppose $s\in S$.  Then $|sN|_S\leq N$.  But by the Claim,
$|sN|_S\geq |s|(|N|_S-4\Delta)$.
Thus $|s|\leq \frac{N}{|N|_S-4\Delta}$. 
\end{proof}
\begin{remark}
Note that the lemma is not true if we do not 
make some assumptions about the geometry of $\Gamma$.
Consider, for example,
 the Cayley graph of $\Z$ with respect to the generating set
$S=\{1,2,4,\ldots,2^n,\ldots\}$.
\end{remark}
\begin{corollary}\label{cor:isometries}
If a group $G$ quasi-acts on a tree $T$ and $g\in G$, then
either $g$ quasi-acts hyperbolically or $g$ quasi-acts elliptically.  
\end{corollary}
\begin{proof}
Define a quasi-action of $\Z$ on $T$ by $n(x)=g^n(x)$.  By Proposition
\ref{prop:quasitree} there is some $S\subset\Z$ and $p\in T$
so that $\Gamma(\Z,S)$
embeds quasi-isometrically and coarsely equivariantly into $T$ via
$n\mapsto g^n(p)$.  By Proposition \ref{prop:Z}, $\Gamma(\Z,S)$ must
either have finite diameter or $S$ must be finite.  If $\Gamma(\Z,S)$
has finite diameter, then so does the orbit $\langle g\rangle p$, 
and so $g$ quasi-acts elliptically.  If $S$ is finite, then
$t\mapsto \lfloor t\rfloor$ is a quasi-geodesic in
$\Gamma(\Z,S)$. As $n\mapsto g^n(p)$ is a quasi-isometric
embedding, it follows that $g$ quasi-acts hyperbolically.
\end{proof}

\subsection{Extracting a pseudocharacter}\label{subsec:pseudo}

\begin{lemma}\label{lemma:almost}
Suppose $\rho\co X\to\R$ is a $(R,\epsilon)$-quasi-isometry, where $X$ is a
graph.  Then there is a $(1,\epsilon')$-quasi-isometry
$\rho'\co X\to\R$ for some $\epsilon'$.
\end{lemma}
\begin{proof}
Suppose that
$\rho\co X\to\R$ is a $(R,\epsilon)$-quasi-isometry.  
By adjusting $\rho$ to be affine on edges and allowing $\epsilon$ to get a
 bit larger, we
may assume that $\rho$ is continuous.
 
As $\rho$ is a
$(R,\epsilon)$-quasi-isometry, the diameter of $\rho^{-1}(0)$ is at most
$\epsilon$.  Furthermore, $X\setminus \rho^{-1}(0)$ has exactly two
unbounded path components, which we denote $P$ and $M$. (We may suppose
that $P\subseteq\rho^{-1}(0,\infty)$ and 
$M\subseteq\rho^{-1}(-\infty,0)$.)  We define a new map
$\rho'\co X\to \R$ by 
\[\rho'(y)=\left\{ \begin{array}{ll}
	d(y,\rho^{-1}(0)) & \mbox{if $y\in P$}\\
	-d(y,\rho^{-1}(0)) & \mbox{if $y\in M$}\\
	0 & \mbox{otherwise}
	\end{array} \right. \]

\begin{claim}\label{claim:4}
The diameter of $\rho'^{-1}(c)$ for $c\in \R$ is at most $4\epsilon$.  The
diameter of $\rho'^{-1}(0)$ is at most $2\epsilon$.
\end{claim}
\begin{proof}
First we bound the diameter of $\rho'^{-1}(0)$.  If $x$ and $y$ are both in
$\rho^{-1}(0)$, then $d(x,y)\leq \epsilon$.  Suppose that
$p\in\rho'^{-1}(0)\setminus\rho^{-1}(0)$.  By continuity, there is some $z\in
P\cup M$ so that $\rho(z)=\rho(p)$. Any path from $z$ to $p$ must of course pass
through $\rho^{-1}(0)$.   Thus 
$d(p,\rho^{-1}(0))<d(z,p)\leq \epsilon$,
and so the diameter of $\rho'^{-1}(0)$ is at most $2\epsilon$. 

Let $c\in\R\setminus\{0\}$, and suppose $x$ and $y$ are in $\rho'^{-1}(c)$.  If
$\rho(x)=\rho(y)$, then $d(x,y)\leq \epsilon$.  Note that $x$ and $y$ are either
both in $P$ or both in $M$, so the signs of $\rho(x)$ and $\rho(y)$ are the
same.  Suppose that $|\rho(y)|>|\rho(x)|$.  Let $z\in \rho^{-1}(0)$, and let
$x'$ be a point on the geodesic from $z$ to $y$ so that $\rho(x')=\rho(x)$. Note
that $d(x,x')\leq \epsilon$.  
Since $z\in \rho^{-1}(0)$ and $\rho'(x)=\rho'(y)$,  $|d(z,x)-d(z,y)|\leq
2\epsilon$. Now,
$d(x,y)\leq d(y,x')+d(x,x')\leq d(y,x')+\epsilon$.  So we're done if we get a
bound on $d(y,x')$.  Note $d(y,z)=d(y,x')+d(x',z)\geq d(y,x')+d(x,z)-\epsilon$.
Rearranging we get $d(y,x')\leq d(y,z)-d(x,z)+\epsilon \leq |d(y,z)-d(x,z)|\leq
3\epsilon$, and so $d(x,y)\leq 4\epsilon$.
\end{proof}

\newcommand{\epsilonprime}{\ensuremath{5\epsilon}}
We now show that
$\rho'$ is a $(1,\epsilonprime)$-quasi-isometry.
Let $x$, $y\in X$.  By the previous claim, we may assume
$\rho'(x)\neq\rho'(y)$.  We then have (up to switching $x$ and $y$) three cases
to consider:  
\begin{enumerate}
\item\label{b} $\rho'(x)=0$
\item\label{c} $\rho'(x)<0<\rho'(y)$
\item\label{d} $\rho'(y)<\rho'(x)<0$ or $0<\rho'(x)<\rho'(y)$ 
\end{enumerate}

In case \ref{b}, $|\rho'(y)-\rho'(x)|=
|\rho'(y)|=d(y,\rho^{-1}(0))$, which differs by at most
$2\epsilon$ from $d(x,y)$.  

In case \ref{c}, note that any path from $x$ to $y$ must pass through 
$\rho^{-1}(0)$.
Let $z\in\rho^{-1}(0)$ lie on some geodesic between $x$ and $y$.
We then have 
\[|\rho'(y)-\rho'(x)|=|\rho'(x)|+|\rho'(y)|\leq d(x,z)+d(y,z)=d(x,y).\]
On the other hand, since $\rho^{-1}(0)$ has diameter at most
$\epsilon$, 
\[|\rho'(y)-\rho'(x)|\geq d(x,z)-\epsilon +d(y,z)-\epsilon = d(x,y)-2\epsilon.\]Thus in this case
$|\rho'(x)-\rho'(y)|$ differs from $d(x,y)$ by at most $2\epsilon$.

In case \ref{d}, let $z\in\rho^{-1}(0)$ be arbitrary.  Then,
\[d(x,y)\geq d(y,z)-d(x,z)\geq |\rho'(y)|-|\rho'(x)|-\epsilon =
|\rho'(y)-\rho'*x)|-2\epsilon.\]
For the remaining inequality, let $x'$ be a point on a geodesic
between $y$ and $z$ so that $\rho'(x')=\rho(x)$.  Then by Claim
\ref{claim:4}, $d(x,x')\leq 4\epsilon$, and so
\begin{eqnarray*}
d(x,y) & \leq & d(x',y)+4\epsilon = d(y,z)-d(x',z)+4\epsilon \\
       & \leq & |\rho'(y)|+\epsilon-|\rho'(x')|+4\epsilon\\
       & =    & |\rho'(y)-\rho'(x)|+5\epsilon
\end{eqnarray*}
\begin{figure}[htbp]
\begin{center}
\begin{picture}(0,0)%
\includegraphics{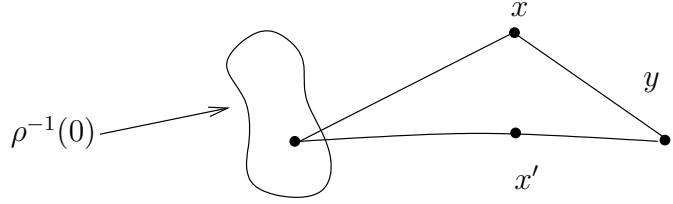}%
\end{picture}%
\setlength{\unitlength}{3947sp}%
\begingroup\makeatletter\ifx\SetFigFont\undefined%
\gdef\SetFigFont#1#2#3#4#5{%
  \reset@font\fontsize{#1}{#2pt}%
  \fontfamily{#3}\fontseries{#4}\fontshape{#5}%
  \selectfont}%
\fi\endgroup%
\begin{picture}(4167,1316)(421,-1255)
\put(4411,-541){\makebox(0,0)[lb]{\smash{{\SetFigFont{12}{14.4}{\familydefault}{\mddefault}{\updefault}{\color[rgb]{0,0,0}$y$}%
}}}}
\put(3578,-98){\makebox(0,0)[lb]{\smash{{\SetFigFont{12}{14.4}{\familydefault}{\mddefault}{\updefault}{\color[rgb]{0,0,0}$x$}%
}}}}
\put(436,-886){\makebox(0,0)[lb]{\smash{{\SetFigFont{12}{14.4}{\familydefault}{\mddefault}{\updefault}{\color[rgb]{0,0,0}$\rho^{-1}(0)$}%
}}}}
\put(3601,-1186){\makebox(0,0)[lb]{\smash{{\SetFigFont{12}{14.4}{\familydefault}{\mddefault}{\updefault}{\color[rgb]{0,0,0}$x'$}%
}}}}
\end{picture}%
\caption{Case \ref{d}: $\rho'(x)$ and $\rho'(y)$ have the same sign.}
\label{fig:sameside}
\end{center}
\end{figure}

Combining the upper and lower
bounds just obtained, we see that in case \ref{d}, $|\rho'(x)-\rho'(y)|$ cannot
differ from $d(x,y)$ by more than $\epsilonprime$.  In all the other cases,
$|\rho'(x)-\rho'(y)|$ is even closer to $d(x,y)$, and so $\rho'$ is a
$(1,\epsilonprime)$-quasi-isometry.
\end{proof}

The remainder of this section is devoted to showing that given a
quasi-action of a group $G$ on a tree $T$, so that one point of
$\partial T$ is fixed by $G$, we may find a
pseudocharacter which is nonzero precisely on those elements of $G$
which act hyperbolically.  In spirit, this pseudocharacter is what we
should expect to get by picking an appropriate Busemann function and
looking at its values on an orbit.  

\begin{proposition}\label{proposition:qchar}
Suppose the finitely generated group $G$ quasi-acts on a simplicial
tree $T$, and suppose that $G$ fixes some point $e\in\partial T$.
Then
there is a pseudocharacter $\chi\co G\to \R$ so that $\chi(g)=0$ if and only
if $g$ quasi-acts elliptically.
\end{proposition}

\begin{proof}
If every $g\in G$ quasi-acts elliptically, we may simply set $\chi=0$.
Otherwise, Corollary \ref{cor:isometries} implies that some 
$\pi\in G$ quasi-acts hyperbolically.
By possibly replacing $\pi$ with $\pi^{-1}$, we may suppose that
$\lim_{n\to\infty} \pi^n x = e$ for any $x\in T$.  Fix $p\in T$ and
let $R$ be large enough so that if $S=\{g\in G\ |\ d(gp,p)\leq R\}$
then $\Gamma=\Gamma(G,S)$ embeds quasi-isometrically in $T$ as in the
conclusion of Proposition \ref{prop:quasitree}.  
By possibly making $R$ a bit larger, we may assume that $\pi\in S$. 
Let $Q\co \Gamma\to
T$ be the quasi-isometric embedding defined by $Q(g)=gp$.  This map is
a quasi-isometry onto its image (a subtree of $T$), and so $\Gamma$ is
a $\delta$-quasi-tree for some $\delta$.  Moreover,
the map $q$ extends to a continuous equivariant
injection $\overline{Q}\co\partial \Gamma
\to\partial T$.  The point $e$ is in the image of $\overline{Q}$,
since
$e=\lim_{n\to\infty}\pi^n(p)=\overline{Q}(\lim_{n\to\infty}\pi^n)$.

Since $\pi\in S$, $\pi^n$ is connected to $\pi^{n+1}$ by an edge for
all $n$ and there is an obvious unit speed  quasi-axis $\gamma\co\R\to \Gamma$ for
$\pi$ satisfying $\gamma(n)=\pi^n$.  This $\gamma$ is a
$(L,D)$-quasi-geodesic for some $L\geq 1$ and $D\geq 0$.

We will begin by finding a ``partial quasicharacter'' on
those group elements which are not too far from
$\gamma$.   We find a sequence of constants
$B_0<B_1<B_2$ which ensure various kinds of good behavior on
$B_i$-neighborhoods of $\gamma$. 
(Note: It will make things slightly easier if we choose the $B_i$ all to be
positive integers, and to satisfy $B_{i+1}>10 B_i$.  This causes no
loss of generality, though of course we will not get the best possible 
constants.)

Since $e$ is fixed, $g\gamma({[0,\infty)})$ is always a finite Hausdorff 
distance from $\gamma({[0,\infty)})$.  The first claim implies that it is
always eventually in some fixed neighborhood of $\gamma$.
\begin{claim}\label{claim:B0}
There is a constant $B_0$ with the following property:  for any $\rho>0$,
there is some $t(\rho)\geq 0$ so that if $d(g,\gamma)\leq \rho$, then
$g\gamma({[t(\rho),\infty)})$ lies in a $B_0$-neighborhood of $\gamma$.
\end{claim}
\begin{proof}(\ref{claim:B0})
The group element $g$ fixes the point $e$ in $\partial \Gamma$
corresponding to $\lim_{t\to\infty}\gamma(t)$.  In particular,
$\lim_{t\to\infty} g\gamma(t)=e$.  
Note that $Q\circ\gamma$ is a continuous
$(KL, KD+C)$-quasi-geodesic in the tree
$T$, and so there is some geodesic $\sigma_1\co\R\to T$ 
whose image is contained in the image
of $Q\circ\gamma$ and which has the same endpoints in 
$\partial T$ as does $Q\circ\gamma$.  
We parameterize $\sigma_1$ so that 
$\lim_{t\to\infty}\sigma_1(t)=\overline(Q)(e)$
tends toward $e$ as $t$ tends to positive infinity.
There is also a geodesic ray $\sigma_2\co{[0,\infty)}\to T$ 
whose image is contained in the image of $Q\circ
g\gamma|_{[0,\infty)}$ and which has the same endpoints in 
$T\cup \partial T$ as does $Q\circ g\gamma|_{[0,\infty)}$.  
Let $\sigma\co{[0,\infty)}\to T$ be the unique geodesic ray whose image
is the intersection of the images of $\sigma_1$ and $\sigma_2$.

Let $t_1 = \inf\{t\ |\ Q\circ\gamma(t)=\sigma(0)\}$ and
$t_g=\inf\{t\ |\ Q\circ g\gamma(t)=\sigma(0)\}$.
Since 
$Q\circ\gamma|_{[t_1,\infty)}$ and $Q\circ g\gamma|_{[t_g,\infty)}$ are
$(KL,KD+C)$-quasi-geodesic rays with the same endpoints as $\sigma$,
neither can stray from a $B:=B(KL,KD+C,0)$-neighborhood of $\sigma$
(Lemma \ref{lemma:uniformnbhd} applied in the case of $T$, a
$0$-hyperbolic space).  Thus for $t\geq t_g$,
the distance between $g\gamma(t)$ and $\gamma$ is at most $K(B+C)$.
(This is because
$Q$ is a $(K,C)$-quasi-isometry and the image of $Q\circ\gamma$
includes all of $\sigma$.)  

Fixing $B_0=KB+KC$, it remains only to show that $t_g$ can be bounded
by some function of $\rho=d(g,\gamma)$.  First note that
$d(Q(g),Q(\gamma))\leq K \rho+C$.  Since $\sigma(0)$ is the closest
point on
$\sigma_1$ to $Q(g)$, and since $Q\circ\gamma$ lies in a
$B$-neighborhood of $\sigma_1$, we have $d(Q(g),\sigma(0))\leq
d(Q(g),Q(\gamma))+B$, and so $d(Q(g),\sigma(0))\leq K\rho+C+B$.
Finally, since 
$d(Q(g),\sigma(0))=d(Q\circ g\gamma(0),Q\circ g\gamma(t_g))$ and
$Q\circ g\gamma$ is a $(KL, KD+C)$-quasi-geodesic, we have 
\[\frac{1}{KL} t_g - (KD+C)\leq d(Q(g),\sigma(0))\]
We may therefore set $t(\rho)=KL(K\rho +C+B+KD+C)$ and the Claim is
established.   
\end{proof}

Let $N_0=\{x\in \Gamma\ |\ d(x,\gamma)\leq B_0\}$, and let $G_0=G\cap
N_0$.  

\begin{claim}\label{claim:B1}
There is some $B_1$ so that if $g\in G_0$, then
$g\gamma|_{[0,\infty)}$ is within $B_1$ of $\gamma$.
\end{claim}
\begin{proof}(\ref{claim:B1})
By Claim \ref{claim:B0}, there is some $t=t(B_0)$ so that if $g$ is
within $B_0$ of $\gamma$, then $g\gamma|_{[t,\infty)}$ is within $B_0$
of $\gamma$.  Thus we need only worry about the initial segment
$g\gamma|_{[0,t]}$.  Since $g\gamma$ is always an
$(L,D)$-quasi-geodesic, it certainly can't go further away than $L
t+D$ in this time, and so any $B_1\geq B_0 + Lt+D$ will suffice. 

(As mentioned before, we want to make sure we have a healthy gap
between $B_0$ and $B_1$, so we choose $B_1>\max\{B_0+Lt+D,10
B_0\}$.)
\end{proof}

Let $N_1=\{x\in \Gamma\ |\ d(x,\gamma)\leq B_1\}$, with the metric induced
by inclusion in $\Gamma$. Notice that this metric may not be the same as
the path metric on $N_1$, so $N_1$ might not be a geodesic metric
space.  We can however find a bigger neighborhood in which geodesics
between different points of $N_1$ may always be found.  

\begin{claim}\label{claim:B2}
There is some $B_2$ so that if $N_2=\{x\in \Gamma\ |\ d(x,\gamma)\leq B_2\}$
is given the \emph{path} metric, then the inclusion of $N_1$ in $N_2$
is an isometric embedding.  That is, between any two points of $N_1$
there is a geodesic in $\Gamma$ between them which lies in a
$B_2$-neighborhood of $\gamma$.   
\end{claim}
\begin{proof}(\ref{claim:B2})
Let $x_1$ and $x_2$ lie in $N_1$, and let $t_1$, $t_2$ be real numbers
so that $d(x_i, \gamma(t_i))\leq B_1$ for $i\in \{1,2\}$.  Since
$\gamma$ is an $(L,D)$-quasi-geodesic, we may use it to define a (not
continuous) $(L, D+2B_1)$-quasi-geodesic $\alpha$ from $x_1$ to $x_2$ 
inside $N_1$ as follows:
\[\alpha(t):=\left\{\begin{array}{ll}
    x_i & t=t_i\\
    \gamma(t) & \mathrm{otherwise}
    \end{array}\right.
\]
Recall that $\Gamma$ is $\delta$-hyperbolic, and so we may make use of
Lemma \ref{lemma:uniformnbhd}; any geodesic in $\Gamma$ between $x_1$
and $x_2$ lies in a
$B(L,D+2B_1,\delta)$-neighborhood of $\alpha$.  Since $\alpha$ lies
entirely inside a $B_1$-neighborhood of $\gamma$, we may set $B_2$ to
some integer at least as large as
$\max\{10 B_1, B_1+B(L,D+2B_1,\delta)\}$.    
\end{proof}

\begin{claim}\label{claim:N2}
$N_2$ is quasi-isometric to $\R$.  
\end{claim}
\begin{proof}(\ref{claim:N2})
By Claim \ref{claim:B2}, any two points in the image of $\gamma$ have
distance in $N_2$ equal to their distance in $\Gamma$.  Thus $\gamma\co
\R\to N_2$ is still an $(L,D)$-quasi-geodesic from the point of view of
the path metric on $N_2$.

Since every point in $N_2$ is within $B_2$ of the image of $\gamma$,
the quasi-isometric embedding $\gamma\co\R\to N_2$ is a quasi-isometry.  
\end{proof}

Since $B_2$ is assumed to be an integer, $N_2$ is a graph, and so
Lemma \ref{lemma:almost} implies that for some $\epsilon >0$ 
there is a continuous
$(1,\epsilon)$-quasi-isometry $\chi_0\co N_2\to \R$.  
By composing with an isometry of $\R$, we may assume that
$\chi_0(1)=0$ and that $\chi_0(\pi^n)\to \infty$ as $n\to \infty$.  
Note that at least
up to some additive error independent of $g$, 
$|\chi_0(g)|$ is the distance from $1$ to $g$ in $\Gamma$.  
We will eventually
show that $\chi_0$ is a ``partial quasicharacter'' on $G_0$
(Claim \ref{claim:partialquasi}).  

If $g\in G_0$, and $\chi_0(hg)-\chi_0(h)$ has the same sign as
$\chi_0(g)$ whenever $hg$ and $h$ are also in $G_0$, we will say $g$
is \emph{irreversible}.  
\begin{claim}\label{claim:irreversible}
There is some $R$ so that if $g\in G_0$ and $|\chi_0(g)|>R$, then $g$ is
irreversible.  
\end{claim}
\begin{proof}(\ref{claim:irreversible})
We prove the claim for $R=(\epsilon+D)L^2+D+B_0$.  
Assume that $g\in G_0$ is not irreversible.  That is, suppose that
there is some $h\in G_0$ so that $gh\in G_0$ but the signs of 
$\chi_0(g)$ and $\chi_0(hg)-\chi_0(h)$ differ.  
Since $g\in G_0$
there is some point $\pi^m$ on $\gamma$ so that $d(g,\pi^m)\leq B_0$.

Suppose the sign of $\chi_0(\pi^m)$ differs from the sign of $\chi_0(g)$.
Since $\pi^m$ and $g$ are both inside $N_1$, a geodesic between
$\pi^m$ and $g$ lies in $N_2$.  By continuity there is some point $w$ on
this geodesic with $\chi_0(w)=0$.  Thus $d(g,1)\leq
d(g,w)+d(w,1)\leq B_0+\epsilon\leq R$.  Thus we may suppose that
the signs of $\chi_0(\pi^m)$ and $\chi_0(g)$ agree.  

Suppose that the sign of $m$ differs from the 
signs of $\chi_0(\pi^m)$ and $\chi_0(g)$.  In this case 
there is some $t\in\R$ having the same sign as $m$, with $|t|>|m|$,
but with $\chi_0(\gamma(t))=0$.  Since $\chi_0(\gamma(t))=0=\chi_0(1)$, 
we must have $d(\gamma(t),1)\leq
\epsilon$, which leads to the following:
\[\frac{1}{L}|t|-D\leq\epsilon
\Rightarrow |t|\leq L(\epsilon+D)
\Rightarrow |m|\leq L(\epsilon+D)\]
\[\Rightarrow d(1,\pi^m)\leq L^2(\epsilon+D)+D
\Rightarrow d(g,1)\leq L^2(\epsilon+D)+D+B_0=R\]
We may therefore suppose that the signs of $m$, $\chi_0(\pi^m)$
and $\chi_0(g)$ all agree.  

We next turn our attention to the order of $\chi_0(h)$, $\chi_0(h\pi^m)$, and
$\chi_0(hg)$.  Suppose first that $\chi_0(h)$ is between $\chi_0(h\pi^m)$
and $\chi_0(hg)$.
Since $h\pi^m$ is at most $B_0$ from $hg$, and
$hg$ is assumed to be at most $B_0$ from $\gamma$, the point $h\pi^m$ is
at most $2B_0<B_1$ from $\gamma$.  
By Claim \ref{claim:B2}, there is a geodesic in $\Gamma$ from $h\pi^m$
to $hg$ entirely inside $N_2$.  
By continuity there is  a point $w$ on
this geodesic with $\chi_0(w)=\chi_0(h)$.  As $\chi_0$ is a
$(1,\epsilon)$-quasi-isometry, $d(w,h)\leq \epsilon$, and so
\[d(g,1)=d(hg,h)\leq d(hg,w)+\epsilon\leq d(hg,h\pi^m)+\epsilon\leq
B_0+\epsilon\]
establishes the claim in this case.  
We may therefore assume that
$\chi_0(h)$ is not between $\chi_0(h\pi^m)$ and $\chi_0(hg)$.

It is now necessary to distinguish cases according to the sign of
$\chi_0(g)$.  We first deal with $\chi_0(g)$ negative and
$\chi_0(hg)-\chi_0(h)$ positive.  Together with our previous
assumptions this gives $\chi_0(h\pi^m)>\chi_0(h)$.  Since
$(h\gamma)^{-1}(h\pi^m)=m$ is negative, $h\gamma$ must
''double back'' as in Figure \ref{fig:positive}
\begin{figure}[htbp]
\begin{center}
\begin{picture}(0,0)%
\includegraphics{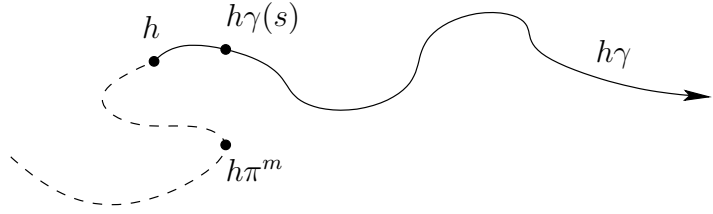}%
\end{picture}%
\setlength{\unitlength}{3947sp}%
\begingroup\makeatletter\ifx\SetFigFont\undefined%
\gdef\SetFigFont#1#2#3#4#5{%
  \reset@font\fontsize{#1}{#2pt}%
  \fontfamily{#3}\fontseries{#4}\fontshape{#5}%
  \selectfont}%
\fi\endgroup%
\begin{picture}(4449,1309)(514,-1424)
\put(1351,-361){\makebox(0,0)[lb]{\smash{{\SetFigFont{12}{14.4}{\familydefault}{\mddefault}{\updefault}{\color[rgb]{0,0,0}$h$}%
}}}}
\put(1876,-1261){\makebox(0,0)[lb]{\smash{{\SetFigFont{12}{14.4}{\familydefault}{\mddefault}{\updefault}{\color[rgb]{0,0,0}$h\pi^m$}%
}}}}
\put(1876,-286){\makebox(0,0)[lb]{\smash{{\SetFigFont{12}{14.4}{\familydefault}{\mddefault}{\updefault}{\color[rgb]{0,0,0}$h\gamma(s)$}%
}}}}
\put(4201,-511){\makebox(0,0)[lb]{\smash{{\SetFigFont{12}{14.4}{\familydefault}{\mddefault}{\updefault}{\color[rgb]{0,0,0}$h\gamma$}%
}}}}
\end{picture}%
\caption{The case of $\chi_0(g)$ negative.  The function $\chi_0$ increases
to the right.  Note that the dotted part
of $h\gamma$ may not lie entirely in $N_2$, though $h\pi^m$ is in $N_2$.}
\label{fig:positive}
\end{center}
\end{figure}
and there must be some $s>0$ with $\chi_0(h\gamma(s))=\chi_0(h\pi^m)$.
Since $\chi_0$ is a $(1,\epsilon)$ quasi-isometry we have 
$d(h\gamma(s),h\pi^m)\leq\epsilon$,
and a by now familiar looking calculation shows that $|s|\leq
L(\epsilon+D)$ and so 
$d(g,1)=d(hg,h)\leq
L^2(\epsilon + D) +D+B_0$.  

In case $\chi_0(g)$ and $m$ are 
positive but $\chi_0(hg)-\chi_0(h)$ is negative, we
again deduce that $h\gamma$ must double back, but now the picture is
slightly different and shown in Figure \ref{fig:negative}.
\begin{figure}[htbp]
\begin{center}
\begin{picture}(0,0)%
\includegraphics{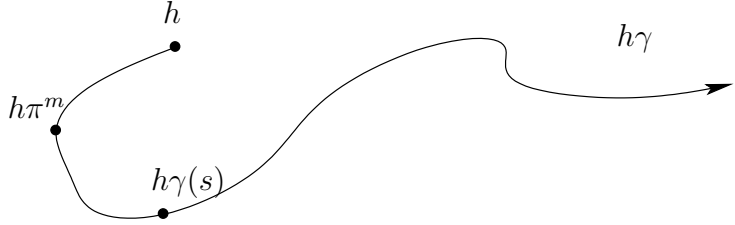}%
\end{picture}%
\setlength{\unitlength}{3947sp}%
\begingroup\makeatletter\ifx\SetFigFont\undefined%
\gdef\SetFigFont#1#2#3#4#5{%
  \reset@font\fontsize{#1}{#2pt}%
  \fontfamily{#3}\fontseries{#4}\fontshape{#5}%
  \selectfont}%
\fi\endgroup%
\begin{picture}(4605,1400)(361,-1602)
\put(1351,-361){\makebox(0,0)[lb]{\smash{{\SetFigFont{12}{14.4}{\familydefault}{\mddefault}{\updefault}{\color[rgb]{0,0,0}$h$}%
}}}}
\put(4201,-511){\makebox(0,0)[lb]{\smash{{\SetFigFont{12}{14.4}{\familydefault}{\mddefault}{\updefault}{\color[rgb]{0,0,0}$h\gamma$}%
}}}}
\put(376,-961){\makebox(0,0)[lb]{\smash{{\SetFigFont{12}{14.4}{\familydefault}{\mddefault}{\updefault}{\color[rgb]{0,0,0}$h\pi^m$}%
}}}}
\put(1276,-1411){\makebox(0,0)[lb]{\smash{{\SetFigFont{12}{14.4}{\familydefault}{\mddefault}{\updefault}{\color[rgb]{0,0,0}$h\gamma(s)$}%
}}}}
\end{picture}%
\caption{The case of $\chi_0(g)$ positive.}
\label{fig:negative}
\end{center}
\end{figure}
Continuity ensures the existence of some $s>m>0$ on $\gamma$ so that
$\chi_0(h\gamma(s))=\chi_0(h)$ (and hence $d(h\gamma(s),h)\leq
\epsilon$).
Now the reader may verify that 
the bound $s-m\leq s \leq L(\epsilon+D)$ leads to 
$d(g,1)\leq L^2(\epsilon+D)+D+B_0$ in this final case.  
\end{proof}

But this is finally enough to prove that $\chi_0$ is a partial
quasicharacter:
\begin{claim}\label{claim:partialquasi}
There is some number $\Delta_0$ so that if $g$, $h$, and $gh$ all lie
in $G_0$, then $|\chi_0(gh)-\chi_0(g)-\chi_0(h)|\leq \Delta_0$.  
\end{claim}
\begin{proof}(\ref{claim:partialquasi})
Suppose $g$, $h$, and $gh$ all lie in $G_0$. 
If $h$ is not irreversible, then $|\chi_0(h)|\leq R$ by Claim
\ref{claim:irreversible}, and so
$d(gh,g)=d(h,1)\leq R+\epsilon$.  Thus since $\chi_0$ is a
$(1,\epsilon)$-quasi-isometry,
$|\chi_0(gh)-\chi_0(g)|=\leq R+2\epsilon$ and
$|\chi_0(gh)-\chi_0(g)-\chi_0(h)|\leq 2R+2\epsilon$ in this case.
  
Suppose that $h$ is irreversible.
Then the sign of $\chi_0(gh)-\chi_0(g)$ is the same as the sign of 
$\chi_0(h)$.  Since
$d(gh,g)=d(h,1)$, the magnitude $|\chi_0(gh)-\chi_0(g)|$ is within 
$2\epsilon$ of
$|\chi_0(h)|$.  Thus $|\chi_0(gh)-\chi_0(g)-\chi_0(h)|\leq 2\epsilon$.  
In both cases
$|\chi_0(gh)-\chi_0(g)-\chi_0(h)|$ is bounded above by $2R+2\epsilon$,
and so $\chi_0$ is a partial 
quasicharacter with defect $\Delta_0\leq 2R+2\epsilon$.
\end{proof}

We now can prove the Proposition by ``extending'' $\chi_0$ to all of
$G$.  For $i>0$, let $G_i=\{g\in G\ |\ \pi^{-i}g\pi^i\in G_0\}$.
We claim that any $g\in G$ is contained in all but finitely many
$G_i$.  Indeed, we have already shown (Claim \ref{claim:B0})
that $g\pi^k$ is always in
$G_0$ for large enough $k$.  Since $\pi G_0 = G_0$, $\pi^{-k}g\pi^k$
is in all $G_i$ for all such $k$.  Thus if we can consistently (at
least up to a bounded error) define partial quasicharacters on all
the $G_i$, we will have defined one on $G$.  Let $\chi_i\co G_i\to \R$
be defined by $\chi_i(g)=\chi_0(\pi^{-i}g\pi^i)$, and define $\chi(g)=
\limsup_{i\to\infty}\chi_i(g)$.  (We could just as well use the
$\liminf$ or any number between the two; the $\liminf$ and $\limsup$
will be seen to be boundedly different from one another, and we only
need understand the coarse properties of $\chi$.)

We now bound the error
$|\delta\chi(g,h)|=|\chi(gh)-\chi(g)-\chi(h)|$ for $g$, $h\in G$.
Note that if $\chi_i(g)$ and $\chi_j(g)$ are both defined for $i<j$, 
then 
\[|\chi_i(g)-\chi_j(g)|=|\chi_0(\pi^{-i}g\pi^i)-
\chi_0((\pi^{i-j})(\pi^{-i}g\pi^i)(\pi^{j-i}))|\leq 3\Delta_0\]
since all powers of $\pi$ also lie in $G_0$.  Thus $\chi(g)-\chi_i(g)$
is also at most $3\Delta_0$, whenever $\chi_i(g)$ is defined.  
Let $k$ be any
number so that $g$, $h$, and $gh$ all lie in $G_k$.  Then
$|\delta\chi(g,h)|\leq |\chi_k(gh)-\chi_k(g)-\chi_k(h)|+9\Delta_0$.
Finally, 
\[|\chi_k(gh)-\chi_k(g)-\chi_k(h)|=|\chi_0(\pi^{-k}gh\pi^k)-\chi_0(\pi^{-k}g\pi^k)-\chi_0(\pi^{-k}h\pi^k)|\leq\Delta_0\]
implies that $|\delta\chi(g,h)|\leq 10\Delta_0$, and so $\chi$ is a
quasicharacter.  As in Remark \ref{remark:pseudoquasi}, if we define
$\overline{\chi}(g)=\lim_{n\to\infty}\frac{\chi(g^n)}{n}$, then
$\overline{\chi}$ is a pseudocharacter which is boundedly different
from $\chi$.  

If $g$ is elliptic, then
$\chi(g^n)$ will be bounded independently of $n$, and so
$\overline{\chi}(g)=0$.

Finally, suppose that $g$ is hyperbolic.  Then either $\{g^n \}$ or $\{g^{-n}
\}$ tend toward point fixed by $G$ in $\partial \Gamma$.  By
replacing $g$ by $g^{-1}$ if necessary, we may assume that the $g^n$ 
tend toward the fixed point as $n\to\infty$.

We now show that the $g^n$ for $n>0$ all lie in some fixed $G_k$.
We construct a quasi-axis for $g$ by
setting
 $\alpha|_{[0,d(1,g)]}$ to
be some geodesic between $1$ and $g$ and requiring that
$\alpha(t+d(1,g))=g\alpha(t)$.  The restriction $\alpha_{[0,\infty)}$
 is a $(L_g,D_g)$-quasi-geodesic
starting at $1$ and limiting on $e$ just as $\gamma_{[0,\infty)}$ is.  
By Lemma \ref{lemma:uniformnbhd} it therefore stays in some
bounded $B_g:=B(\max\{L,L_g\}, \max\{D,D_g\},\delta)$-neighborhood of 
$\gamma$.
Together with Claim \ref{claim:B0}, this
implies that there is some $N$ so that $g^n\pi^k\in G_0$ for all
$n>0$, $k>N$.  Thus all positive powers of $g$ lie in $G_N$ for some
fixed $N$ which depends only on $g$.  Since
$g$ acts hyperbolically, $d(1,\pi^{-N}g^n\pi^{N})$ tends to infinity
as $n$ tends to infinity, implying that $\chi_N$ is unbounded on 
$\{g^n\ | \ n>0\}$. 
This implies  that $\chi$ and $\overline{\chi}$ are
likewise unbounded on $\{g^n\ |\ n>0\}$.  Since $\overline{\chi}$ is a
homomorphism on cyclic subgroups, we must
have $\overline{\chi}(g)\neq 0$. 
\end{proof}

\section{A class of groups with Property \QFA }\label{section:main}

In this section we give an example of how the results of Section
\ref{section:lemmas} can be used to establish Property \QFA.

\begin{definition}
Let $G$ be a group, and let $g$ be an element of $G$.  We will say $g$
is a \emph{stubborn element of $G$} if for all
$H< G$ with $[G:H]\leq 2$, there exists some integer $k_H>0$ so that
$g^{k_H}\in [H,H]$.  
\end{definition}

\begin{lemma}\label{lemma:Bgeneral}
Let $B$ be an amenable group, quasi-acting on a tree $T$.
If $b$ is a stubborn element of $B$, then $b$ quasi-acts
elliptically.
\end{lemma}
\begin{proof}
Suppose that $b$ does not quasi-act elliptically.
We show that $B$ cannot be amenable.
By Corollary \ref{cor:isometries} $b$ quasi-acts hyperbolically, and
therefore fixes exactly two points in $\partial T$.  Let
$\mathrm{Fix}(b)=\{e_1,e_2\}$ be the fixed point set of $b$.

We claim that either $B$ preserves $\mathrm{Fix}(b)$ or fixes either
$e_1$ or $e_2$.  Indeed,
suppose that $\mathrm{Fix}(b)$ is not preserved by all of $B$.  If
there are elements $b_1$ and $b_2$ so that $b_1(e_1)\neq e_1$ and
$b_2(e_2)\neq  e_2$, then 
a simple ping-pong argument shows that $B$ must contain a nonabelian
free subgroup and thus cannot be amenable.

We may therefore assume that there is some subgroup $H\leq B$ of index
at most $2$ which 
fixes some point in $\partial X$ and contains the element $b$.
By Proposition
\ref{proposition:qchar} there is a pseudocharacter $p\co H\to \R$ so that
$p(b)\neq 0$.  Since $b^{k_H}\in[H,H]$, $p$ is not a homomorphism, and so 
$[\delta p]$ is a nonzero element of $H_b^2(H,\R)$.  Theorem
\ref{theorem:trauber} then implies that $H$ is not amenable, and so
$B$ is not amenable.  
\end{proof}
\begin{definition}
A group $G$ is \emph{boundedly generated} if there exists a finite tuple
$(g_1,\ldots,g_n)$  of elements of $G$ so that any $g\in G$ is equal to
$g_1^{\alpha_1}\cdots g_n^{\alpha_n}$ for some tuple
$(\alpha_1,\ldots,\alpha_n)$ of integers.  This is easily seen to be
equivalent to the following
condition:  There is some finite set $\{g_1,\ldots,g_n\}\subset G$ 
so that the Cayley
graph $\Gamma(G,S)$ has finite diameter
for any
$S\supseteq\{g_i^m\ |\ i\in\{1,\ldots,n\},\ m\in\Z\}$.  
In either case we say that
$G$ is boundedly generated by the elements $g_1\ldots,g_n$.  
\end{definition}

\begin{theorem}\label{theorem:qfa}
Let $G$ be a group which is boundedly generated by elements 
$g_1,\ldots,g_n$, so
that for each $i$, $g_i$ is a stubborn element of $B_i$ for some
$B_i < G$. Then $G$ has Property \QFA.
\end{theorem}
\begin{proof}
Suppose that $G$ quasi-acts on some tree $T$.  
By Lemma \ref{lemma:Bgeneral}, each $g_i$ quasi-acts elliptically.

Fix $p\in T$. We must show that the orbit $Gp$ is bounded.
Let 
\[R > \sup\{d(p,g_i^m p)\ |\ {i\in \{1,\ldots,n\},\  m\in \Z}\}.\]  Since
each $g_i$ has bounded orbits, we may find such a finite $R$.  
By perhaps increasing $R$, we may assume that $R$ is large enough that
$\Gamma(G,S)$ coarsely equivariantly quasi-isometrically embeds in $T$, for
$S=\{g\in G\ |\ d(gp,p)\leq  R\}$ (applying Proposition \ref{prop:quasitree}).  
But since $S\supseteq\{g_i^m\ |\ i\in\{1,\ldots,n\},\ m\in\Z\}$, the Cayley
graph $\Gamma(G,S)$ has finite diameter.  It follows that the orbit $Gp$ has
finite diameter.
\end{proof}

\begin{corollary}
Let $K$ be a number field, and $\mathcal{O}$ its ring of integers.  Then
$SL(n,\mathcal{O})$ has Property \QFA, for $n>2$.  In particular, $SL(n,\Z)$ has
property \QFA, for $n>2$.
\end{corollary}
\begin{proof}
Let $\Lambda=\{\lambda_1,\ldots,\lambda_k\}$ be an integral basis for
$\mathcal{O}$ (for $\mathcal{O}=\Z$, $\Lambda=\{1\}$).  
Let $e_{ij}^\lambda$ be the matrix equal to the identity
matrix except for the $ij$th entry, which is equal to $\lambda$.  
Note that
$(e_{ij}^\lambda)^{-1}=e_{ij}^{-\lambda}$.
It is shown in
\cite{carterkeller:slno} 
that $SL(n,\mathcal{O})$ is boundedly generated by 
$\{e_{ij}^\lambda\ |\ i\neq j,\ \lambda\in\Lambda\}$.
It can be easily verified that $e_{ik}^\lambda e_{kj}^1 e_{ik}^{-\lambda}
e_{kj}^{-1}=e_{ij}^\lambda$, provided that $i$, $j$ and $k$ are all distinct. 
Furthermore, the subgroup $B_{ij}=\langle e_{ik}^\lambda, e_{kj}^1\rangle$ is
nilpotent and therefore amenable.  If $H<B_{ij}$ is a subgroup of
index $2$, note that $e_{ik}^{2\lambda}$ and $e_{kj}^{2}$ are
contained in $H$.  The commutator $[e_{ik}^{2\lambda},e_{kj}^2] =
e_{ij}^{4\lambda} = (e_{ij}^{\lambda})^4$, so $e_{ij}^{\lambda}$ is stubborn.
Theorem \ref{theorem:qfa} can then be applied.
\end{proof}

\begin{remark}
It is natural to wonder to what extent the hypotheses of Theorem
\ref{theorem:qfa} can be weakened or modified, and what other groups
can be shown to have Property \QFA\ using these methods.  
Using Theorem \ref{theorem:qfa}
together with
theorems of Tits \cite{tits:congruence} and Tavgen$'$
\cite{tavgen:bg}
it is
possible to prove Property \QFA\ holds for simple, simply connected
Chevalley groups of semisimple rank at least two
over the ring of integers of a number field.

One might wonder 
whether bounded generation plus \FA\ is enough to guarantee \QFA.
In fact there are groups
(certain central extensions of lattices in higher rank Lie groups) which
are boundedly generated and have Property \FA\ (in fact they are
Kazhdan, which is stronger), but which
admit a nontrivial pseudocharacter.  Thus by the construction in 
Example \ref{example:quasifrompseudo} these groups
do not have Property \QFA.  These examples are described in the
Appendix.
\end{remark}

\section{\QFA-type properties for particular kinds of trees}\label{section:coda}

We can break Question \ref{question:main} down into parts by considering only
cobounded actions on particular types of trees.  That is the strategy in this
section.    
\begin{definition}
If $\mathcal{T}$ is a class of (infinite diameter)
trees, and $G$ does not admit a cocompact action on any tree in $\mathcal{T}$,
then we will say that $G$ has Property \fa{T}.  If $G$ does not admit a
cobounded quasi-action on any tree in $\mathcal{T}$, we will say that $G$ has
Property \qfa{T}.
\end{definition}

Note that while \qfa{T}\ always implies \fa{T}, \fa{T}\ does not usually imply
\qfa{T}.

\begin{definition}
A tree $T$ is \emph{bushy} if there is a number $B$, called the 
\emph{bushiness 
constant} of $T$, so that for any point $p\in T$ the set $\{x\in T\ |\
d(p,x)\geq B\}$ has at least three components.  We will say that $T$ is
\emph{finitely bushy} if in addition $\{x\in T\ |\ d(p,x)\geq B\}$
 always has finitely many components.  If $B$ can be chosen so that 
$\{x\in T\ |\ d(p,x)\geq B\}$
has a countable infinity of components for every $x\in T$,
we say the tree is
\emph{countably bushy}.

We write $\mathcal{B}$ for the set of 
finitely bushy trees and $\mathcal{C}$ for 
the set of countably bushy trees.
\end{definition}

The proof of the following proposition is left to the reader:
\begin{proposition}
A finitely generated 
group $G$ has Property \QFA\ if and only if it has all three of the properties
\qfar, \qfa{B}\ and \qfa{C}.  
\end{proposition}

\begin{example}
The group $\Z$ obviously does not have \qfar, but Propositions
\ref{prop:quasitree}
and \ref{prop:Z} imply that $\Z$ has \qfa{B}\ and \qfa{C}.
\end{example}

It is a corollary of the main theorem of \cite{msw:quasiactI} that Properties
\fa{B}\ and \qfa{B}\ are equivalent.  A characterization of \qfar\ is given by
Proposition \ref{prop:qfar} below.  
\begin{proposition}\label{prop:qfar}
Let $G$ be a group.  The following are equivalent:
\begin{enumerate}
\item\label{qfar} $G$ has Property \qfar.
\item\label{pchar} Neither $G$ nor any index $2$ subgroup of $G$ has a 
nontrivial pseudocharacter.
\end{enumerate}
\end{proposition}
\begin{proof}

\textbf{\ref{qfar} implies \ref{pchar}:}

We argue by assuming 
 that \ref{pchar} does not hold and then producing a quasi-action by $G$ on
$\R$. 

Suppose either $G$ or some index two subgroup of $G$ admits a nontrivial pseudocharacter.  The easiest case to deal with is if $f\co G\to\R$ is a nontrivial
pseudocharacter.  Then $A\co G\times \R\to\R$ defined by $A(g,x)=f(g)+x$ is a
cobounded quasi-action on $\R$.
 
Let $H$ be an index $2$ subgroup of $G$ and suppose that $f\co H\to\R$ is a
nontrivial pseudocharacter.
Fix some element $t$ of $G\setminus H$ and define $f_t\co H\to G$ by
$f_t(h)=f(h)-f(t^{-1}h t)$.  Note that $f_t(t^{-1}ht)=-f_t(h)$.  
There are two cases, depending on whether $f_t(h)=0$ for every $h$ in $H$.  

If $f_t=0$ on $H$ we can extend $f$ to a quasicharacter
$\bar{f}\co G\to \R$ by setting:
\[\bar{f}(g)=\left\{\begin{array}{ll}
	f(g)     &      \mbox{if $g\in H$}\\
	f(h)     &      \mbox{if $g=th$ for $h\in H$}
	\end{array}
\right.\]
To see that $\bar{f}$ is a quasicharacter, we must bound
 $|\delta \bar{f}(\alpha,\beta)|$
in four cases, depending on whether $\alpha$ and $\beta$ are in $H$ or $tH$.
In three cases, the reader may verify that $|\delta\bar{f}(\alpha,\beta)|$ is
bounded by $\|\delta f\|$.  In the fourth case, where $\alpha$ and $\beta$ both
lie in $tH$,  suppose that $\alpha=th_1$ and $\beta=th_2$.  Then
$\bar{f}(\alpha\beta)=\bar{f}(th_1th_2)=\bar{f}(t^2(t^{-1}h_1t)h_2)$, which
differs from $\bar{f}(\alpha)+\bar{f}\beta$ by at most $2\|\delta
f\|+|f(t^2)|$.  Thus we have produced an unbounded quasicharacter on $G$, which
may be modified to a pseudocharacter as in Remark \ref{remark:pseudoquasi}, and
we are back in the easy case.
 
If $f_t$ is not identically zero, then it is a nontrivial pseudocharacter on
$H$.  It is then possible to build a quasi-action of $G$ on $\R$ so that $t$
acts as a reflection.  Namely, we set
\[A(g,x)=\left\{\begin{array}{ll}
	f_t(g)+x     &     \mbox{if $g\in H$}\\
	-f_t(h)-x    &     \mbox{if $g=th$ for $h\in H$}
	\end{array}
\right.\]
Each group element acts as an isometry, so we only need to show that
$(\alpha)(\beta x)$ and $(\alpha\beta)x$ are uniformly close.  
Again there are
four different cases, depending on whether $\alpha$ and $\beta$ lie in $H$ or
$tH$.  All four cases may safely be left to the diligent reader.  Since $f$ is
assumed to be nontrivial, it is unbounded on $H$, and thus this quasi-action is
cobounded.

\textbf{\ref{pchar} implies \ref{qfar}:}

We assume that there is a cobounded quasi-action of $G$ on $\R$ and produce a
pseudocharacter.

Any element of $G$ must either switch $\pm\infty$ or preserve them.  Thus either
$G$ or an index two subgroup of $G$ fixes 
 $\pm\infty$.  We may then apply
Proposition \ref{proposition:qchar} to obtain a nontrivial
pseudocharacter on the subgroup of $G$
fixing $\pm\infty$.  
\end{proof}

Of course it remains to give a satisfactory account of \qfa{C}.

\def\cprime{$'$}

\newpage
\appendix
\section{Boundedly generated groups with pseudocharacter(s) \protect{\\} by
  Nicolas Monod and Bertrand R\'emy}
\def\R{\mathbb{R}}
\def\Z{\mathbb{Z}}
\def\Q{\mathbb{Q}}
\def\h{\mathrm{H}}
\def\hb{\mathrm{H}_\mathrm{b}}
\def\hc{\mathrm{H}_\mathrm{c}}
\def\hbc{\mathrm{H}_\mathrm{cb}}

\noindent
The aim of this appendix is to construct concrete groups which simultaneously:

\smallskip

(1)~are boundedly generated;

(2)~have Kazhdan's property~(T);

(3)~have a one-dimensional space of pseudocharacters.

\smallskip

By~(3), such groups don't have property~\QFA, whilst they have
property~\FA\ by~ (2);
moreover the quasimorphisms in~(3) cannot be \emph{bushy} in the
sense of~\cite{Manning1}.
Property~(3) has its own interest, as all previous constructions
yield infinite-dimensional spaces. (By taking direct products of our
examples, one gets any finite
dimension.) The examples will be lattices $\widetilde{\Gamma}$ in
non-linear simple Lie groups; more
precisely, starting with certain higher rank Lie groups $H$ with
$\pi_1(H) = \Z$ and suitable lattices
$\Gamma<H$, the group $\widetilde{\Gamma}$ will be the preimage of
$\Gamma$ in the universal covering
central extension
$$0\longrightarrow \Z \longrightarrow \widetilde{H} \longrightarrow H
\longrightarrow 1.\eqno{(*)}$$

\bigskip

Let us first start with \emph{any} group $\Gamma$ satisfying the
following cohomological
properties (we refer to~\cite{Burger-Monod3} for our use of bounded
cohomology):

\smallskip

(3$'$)~the second bounded cohomology $\hb^2(\Gamma,\R)$ has dimension one;

(3$''$)~the natural map
$\psi_\Gamma:\hb^2(\Gamma,\R)\to\h^2(\Gamma,\R)$ is injective;

(3$'''$)~the image of the natural map
$i_\Gamma:\h^2(\Gamma,\Z)\to\h^2(\Gamma,\R)$ spans the image of $\psi_\Gamma$.

\smallskip

We claim that under these assumptions, there is a central extension
$0 \to \Z \to \widetilde\Gamma \to \Gamma \to 1$
such that the kernel of $ \psi_{\widetilde\Gamma} :
\hb^2(\widetilde\Gamma,\R) \to \h^2(\widetilde\Gamma,\R)$
has dimension one.

\smallskip

{\it Proof.~}
By the assumptions, there is $\omega_\Z \!\in\! \h^2(\Gamma,\Z)$
and $\omega \!\in\! \hb^2(\Gamma,\R)$ such that
$\psi_\Gamma(\omega)=i_\Gamma(\omega_\Z)\neq 0$. The central extension
$0 \longrightarrow \Z \longrightarrow \widetilde\Gamma {\buildrel\pi
\over \longrightarrow} \Gamma \longrightarrow 1$ associated to
$\omega_\Z$ yields a commutative diagram:

\smallskip

\centerline{$\xymatrix{
\hb^2(\Gamma,\R) \ar[r]^{\psi_\Gamma} \ar[d]^{\pi^{*}_{\mathrm{b},\R}}
& \h^2(\Gamma,\R) \ar[d]^{\pi^{*}_\R} & \h^2(\Gamma,\Z)
\ar[d]^{\pi^{*}_\Z} \ar[l]_{i_\Gamma}
\\
\hb^2(\widetilde\Gamma,\R) \ar[r]^{\psi_{\widetilde\Gamma}} &
\h^2(\widetilde\Gamma,\R) &
\h^2(\widetilde\Gamma,\Z)
\ar[l]_{i_{\widetilde\Gamma}}}$}

\smallskip

Since $\Z$ is amenable, $\pi^*_{\mathrm{b},\R}$ is an
isomorphism~\cite[3.8.4]{Ivanov}
(this is not true in general for $\Z$ coefficients). Setting
$\beta:=\pi^{*}_{b,\R}(\omega)$, we are reduced to
seeing that $\hb^2(\widetilde\Gamma,\R) = \R \beta$ maps trivially to
$\h^2(\widetilde\Gamma,\R)$.
But we have:
$\psi_{\widetilde\Gamma}(\beta)
=\pi^{*}_{\R} \bigl( \psi_\Gamma(\omega) \bigr)
=(\pi^{*}_{\R} \circ i_\Gamma)(\omega_\Z)
=(i_{\widetilde\Gamma} \circ \pi^{*}_{\Z})(\omega_\Z)$, and
$\widetilde\Gamma$ was designed as a central extension in order to have
$\pi^{*}_{\Z}(\omega_\Z)=0$.
\proofbox\smallskip

{\it Remarks.~}
1. $\widetilde\Gamma$ has property~(T) whenever $\Gamma$ does. Indeed,
since $\psi_\Gamma(\omega)\neq 0$, we have $\omega_\Z\neq 0$ and the
corresponding central extension does not split. The claim is now
a result due to Serre~\cite[p.~41]{Harpe-Valette}.

2. The space of pseudocharacters of $\widetilde\Gamma$ is isomorphic
to $\mathrm{Ker}(\psi_{\widetilde\Gamma})$
modulo the characters of $\widetilde\Gamma$; in particular, since
property~(T) groups have no non-zero characters, $\widetilde\Gamma$
satisfies~(3)
if $\Gamma$ was chosen with property~(T).

3. The group $\widetilde\Gamma$ is boundedly generated whenever $\Gamma$ is so.

\medskip

In conclusion, it remains to check the existence of groups $\Gamma$
satisfying (1),
(2) and (3$'$)-(3$'''$).
We obtain two families of examples from the following discussion
(see also Remark 4 below).

\smallskip

Let $X$ an irreducible Hermitian symmetric space of non-compact type.
Let $H:={\rm Isom}(X)^\circ$ be the identity component of its
isometry group.
We assume that  $\pi_1(H) = \Z$, i.e. that $\pi_1(H)$ is torsion-free.
We have then a central
extension as in~$(*)$ above, yielding a class $\omega_{H,\Z}$ in
the ``continuous'' cohomology $\hc^2(H,\Z)$ (represented by a Borel cocycle);
the image $\omega_H$ of $\omega_{H, \Z}$ under the natural map
$\hc^2(H,\Z)\to \hc^2(H,\R)$ generates $\hc^2(H,\R)$. For all
this, see~\cite{Guichardet-Wigner}.

Let now $\Gamma<H$ be any lattice and let $\omega_\Z$ be the image
of $\omega_{H, \Z}$ under the restriction map
$r_\Z:\hc^2(H,\Z)\to \h^2(\Gamma,\Z)$; thus, the
corresponding central extension
$\widetilde\Gamma$ is (isomorphic to) the preimage of $\Gamma$ in
$\widetilde{H}$.
Note that, so far, $\omega_\Z$ can be zero. From now on we
assume that the rank of $X$ is at least two.
This implies on one hand that $H$ and $\Gamma$ have property~(T)~\cite[2b.8 and
3a.4]{Harpe-Valette}; on the other hand,
  (3$''$) is established in~\cite[Thm.~21]{Burger-Monod3}.
Furthermore, there are isomorphisms
$\hc^2(H,\R){\buildrel\psi\over\longleftarrow}\hbc^2(H,\R){\buildrel{r_\R}\over\longrightarrow}\hb^2(\Gamma,\R)$
(see~\cite{Burger-Monod3} for
the first and the vanishing theorem in~\cite{Monod-Shalom1} for the second).
Thus (3$'$) and (3$'''$) follow as well given the above discussion of
the cohomology of $H$.

\smallskip

Finally, we investigate when $\Gamma$ (and
thus $\widetilde\Gamma$) can be chosen to be boundedly generated
using a result of
Tavgen'~\cite[Theorem B]{Tavgen}.
We define $\Gamma$ as integral points of a $\Q$-algebraic group
$\underline H$ such that the
identity component $\underline H(\R)^\circ$ is $H={\rm Isom}(X)^\circ$.
Using Tavgen's theorem requires that $\underline H$ be quasi-split over $\Q$.
According to Cartan's classification~\cite[X.6, Table V and \S
3]{Helgason78}, the
exceptional cases $E\,III$ and $E\,VII$, and the classical series $D\,III$, are
excluded because the isometry groups are not quasi-split, and {\it a
fortiori~} neither
are their $\Q$-forms.
Let us check that the remaining types admit quasi-split $\Q$-forms.

\smallskip

Case $C\, I$: this corresponds to Siegel's upper half-spaces ${\rm
Sp}_{2n}(\R)/{\rm U}(n)$.
The standard symplectic forms with all coefficients equal to one
define $\Q$-split
algebraic subgroups of ${\rm SL}_{2n}$~\cite[V.23.3]{Borel91}.
For each $n \geq 2$, the lattice
$\Gamma={\rm Sp}_{2n}(\Z):={\rm Sp}_{2n}(\Q) \cap {\rm SL}_{2n}(\Z)$
satisfies all the
required properties. The corresponding symmetric space $X$ has rank $n$ and
dimension $n(n+1)$.

Case $A\, III$: this corresponds to
${\rm SU}(p,q)/{\rm S}\bigl({\rm U}(p)\times{\rm U}(q)\bigr)$ with $p \geq q$.
In view of the Satake-Tits diagrams~\cite[II \S 3]{Satake}, the corresponding
isometry groups which are quasi-split over $\R$ are those for which
$p=q$ or $p=q+1$.
The Hermitian form $h:=\bar x_1x_{2n} - \bar x_2x_{2n-1} + ... -
x_1\bar x_{2n}$
(resp. $\bar x_1x_{2n+1} - \bar x_2x_{2n} + ... - x_1\bar x_{2n+1}$), where
$\bar{}$ denotes the conjugation of $\Q(i)$, defines a $\Q$-form of
the isometry
group ${\rm SU}(n,n)$ (resp. ${\rm SU}(n+1,n)$).
The matrices of ${\rm SL}_{2n}(\Z[i])$ (resp. ${\rm
SL}_{2n+1}(\Z[i])$) preserving $h$ provide
suitable groups $\Gamma$.

\smallskip

{\it Remarks.~}
1. What we call {\it bounded generation}, following \emph{e.g.}~\cite[\S A.2
p.575]{Platonov-Rapinchuk} and
\cite{ShalomIHES}, is what Tavgen' calls \emph{finite width}, while
bounded generation in~\cite{Tavgen}
is defined with respect to a generating system.

2. To have bounded generation, we restricted ourselves to arithmetic
subgroups of quasi-split groups,
which prevents from constructing the groups $\Gamma$ as uniform
lattices (the Godement compactness criterion
requires $\Q$-anisotropic 
groups~\cite[Theorem~4.12]{Platonov-Rapinchuk}, which are so to
speak opposite to split and quasi-split groups).
The underlying deeper problem is to know whether boundedly
generated \emph{uniform} lattices exist~\cite[Introduction]{ShalomIHES}.

3. Given the cohomological vanishing results of~\cite{Burger-Monod3}, 
\cite{Monod-Shalom1},
the only possibilities for $\Gamma$ to be a lattice in (the 
$k$-points of) a simple group over
a local field $k$ is the case we considered: $k=\R$, rank at least 
two and Hermitian structure.
In particular, the non-Archimedean case is excluded. As far as 
bounded generation only is
concerned, there is an even stronger obstruction in positive 
characteristic: any boundedly
generated group that is linear in positive characteristic is virtually
Abelian~\cite{Abert-Lubotzky-Pyber}.

4. A case in Cartan's classification was not alluded to above.
This is the type $BD\, I$, corresponding to
${\rm SO}(p,q)^\circ/\bigl({\rm SO}(p)\times{\rm SO}(q)\bigr)$  with
$p \geq q=2$.
First, ${\rm SO}(2,2)^\circ$ is not simple and the associated symmetric space
is not irreducible (it is the product of two hyperbolic disks).
For $p \geq 3$, the fundamental group $\pi_1({\rm SO}(p,2)^\circ)$ has torsion
since it is $\Z \oplus \Z/2\Z$ ~\cite[I.7.12.3]{Husemoller},
but lattices in $H={\rm SO}(p,2)^\circ$ still enjoy properties (2) and
(3$'$)-(3$'''$).
For bounded generation, since a symmetric non-degenerate bilinear 
form defines a split (resp.
quasi-split) orthogonal group if and only if
$p-q \leq 1$ (resp. $p-q \leq 2$)~\cite[V.23.4]{Borel91}, suitable 
groups $\Gamma$ are
provided by lattices ${\rm SO}(Q) \cap {\rm SL}_n(\Z)$, with $Q$
a non-degenerate quadratic form on $\Q^n$ of signature
$(3,2)$ or $(4,2)$ over $\Q$.

\def\cprime{$'$}
{\renewcommand{\refname}{Appendicular references}

\end{document}